\documentclass[11pt]{amsart}
\usepackage{amssymb,amsmath,amsthm}
\usepackage{mathrsfs,a4wide}
\theoremstyle{plain}
\newtheorem{theorem}{Theorem}[section]
\newtheorem{lemma}[theorem]{Lemma}
\newtheorem{proposition}[theorem]{Proposition}

\newtheorem{remark}[theorem]{Remark}

\theoremstyle{definition}
\theoremstyle{remark}
\numberwithin{equation}{section}

\newcommand{\as}{{\mathcal A}}
\newcommand{\hs}{{\mathcal H}}
\newcommand{\ks}{{\mathcal K}}
\newcommand{\cs}{{\mathcal C}}

\newcommand{\fs}{{\mathcal F}}

\newcommand{\ns}{{\mathcal N}}
\newcommand{\bs}{{\mathcal B}}

\newcommand{\Es}{{\mathcal E}}
\newcommand{\ts}{{\mathcal T}}
\newcommand{\rs}{{\mathcal R}}

\newcommand{\tb}{{\bf T}}
\newcommand{\rb}{{\bf R}}
\newcommand{\seub}{{\bf S}}

\newcommand{\R}{{\mathbb R}}
\newcommand{\N}{{\mathbb N}}

\newcommand{\tint}[1]{{\rm int}(#1)}

\newcommand{\Om}{\Omega}
\newcommand{\Omb}{\overline{\Omega}}

\newcommand{\weakst}{\stackrel{\ast}{\rightharpoonup}}

\newcommand{\weak}{\rightharpoonup}

\newcommand{\eps}{\varepsilon}
\newcommand{\afeaom}{\as_{\varepsilon,a}(\Om)}
\newcommand{\afenaom}{\as_{\varepsilon_n,a}(\Om)}
\newcommand{\afeom}{\as \fs_{\varepsilon}(\Om)}
\newcommand{\afenom}{\as \fs_{\varepsilon_n}(\Om)}

\newcommand{\treaom}{\ts_{\varepsilon,a}(\Om)}

\newcommand{\Sg}[2]{S^{#1}(#2)}
\newcommand{\vuat}{u_a(t)}
\newcommand{\vuas}{u_a(s)}

\newcommand{\zr}{z^r}

\newcommand{\ur}{u^r}

\newcommand{\unr}{u^r_n}

\newcommand{\znr}{z^r_n}

\newcommand{\ulp}{(z^l)^+}
\newcommand{\ulm}{(z^l)^-}
\newcommand{\umk}{z^{r(k)}}
\newcommand{\umkp}{(z^{r(k)})^+}
\newcommand{\umkm}{(z^{r(k)})^-}
\newcommand{\umkn}{z^{r(k)}_n}

\newcommand{\cint}{K}
\newcommand{\Enk}{E_n^{k}}
\newcommand{\Enkp}{E^{k,+}_n}
\newcommand{\enkpt}{\tilde{E}^{k,+}_n}
\newcommand{\enkmt}{\tilde{E}^{k,-}_n}
\newcommand{\Enkm}{E^{k,-}_n}
\newcommand{\rskn}{\rs_n(Q_k)}
\newcommand{\hkn}{H_{k}}
\newcommand{\hknp}{H^{k,+}_n}
\newcommand{\hknm}{H^{k,-}_n}

\newcommand{\cks}{c_s(r(k))}
\newcommand{\cksu}{c_1(r(k))}
\newcommand{\cksd}{c_2(r(k))}

\newcommand{\wsig}{w_\sigma}
\newcommand{\twn}{\tilde{w}_n}
\newcommand{\hwn}{\hat{w}_n}

\input psfig.sty

\newcommand{\res}{\mathop{\hbox{\vrule height 7pt width .5pt depth 0pt
\vrule height .5pt width 6pt depth 0pt}}\nolimits}

\title
[A discontinuous finite element approximation
]
{A discontinuous finite element approximation of \\
quasi-static growth of brittle fractures}
\author[A. Giacomini]
{Alessandro Giacomini}
\address[Alessandro Giacomini]{S.I.S.S.A., Via Beirut 2-4, 34014, Trieste,
Italy}
\email[A. Giacomini]{giacomin@sissa.it}
\author[M. Ponsiglione]
{Marcello Ponsiglione}
\address[Marcello Ponsiglione]{S.I.S.S.A., Via Beirut 2-4, 34014, Trieste,
Italy}
\email[M. Ponsiglione]{ponsigli@sissa.it}
\begin{document}
\vskip .2truecm
\begin{abstract}
\small{
We propose a discontinuous finite element approximation for a model
of quasi-static growth of brittle fractures in linearly elastic bodies formulated
by Francfort and Marigo, and based on the classical Griffith's criterion. We restrict our
analysis to the case of anti-planar shear and we consider discontinuous displacements
which are piecewise affine with respect to a regular triangulation.
}
\end{abstract}
\maketitle
{\small \tableofcontents}

\section{Introduction}
\label{intr}

In this paper we formulate a discontinuous finite element approximation
for a model of quasi-static growth of brittle fractures in linearly
elastic bodies proposed by Francfort and Marigo \cite{FM}. Their
model is based on the classical Griffith's criterion which involves a
competition between {\it bulk} and {\it surface} energies. To be precise,
let $\Om \subseteq \R^3$ be an elastic body, $\partial_D \Om$ be a part of its
boundary and let $g: \partial_D \Om \to \R^3$ be the spatial displacement of $\Om$ at
the points of $\partial_D \Om$. According to Griffith's theory,
given a preexisting crack $\Gamma_1 \subseteq \Omb$,
the new crack $\Gamma$ and the displacement $u: \Om \setminus \Gamma \to \R^3$
associated to $g$ at the equilibrium minimizes the following
energy
\begin{equation}
\label{griffithenergy}
\Es(v,g,\Gamma):= \int_\Om \mu |Ev|^2 +\frac{\lambda}{2}\,|{\rm div}\,v|^2 \,dx+ \hs^2(\Gamma),
\end{equation}
among all cracks $\Gamma$ with $\Gamma_1 \subseteq \Gamma$ and all displacements
$v: \Om \setminus \Gamma \to \R^3$ with $v=g$ on $\partial_D \Om \setminus \Gamma$. Here $Ev$
denotes the symmetric part of the gradient of $v$, and $\hs^2$ denotes the two
dimensional Hausdorff measure, while $\mu$ and $\lambda$ are the Lam\'e coefficients.
The boundary condition is required only on $\partial_D \Om \setminus \Gamma$
because the displacement in a fractured region is assumed not to be transmitted.
Let us indicate by $\Es(g,\Gamma)$ the minimum value of
(\ref{griffithenergy}) among all $v: \Om \setminus \Gamma \to \R^3$
with $v=g$ on $\partial_D \Om \setminus \Gamma$.
\par
Supposing that the boundary displacement $g$ varies with the time $t \in [0,1]$, the quasi-static
evolution $t \to \Gamma(t)$ proposed in \cite{FM} requires that:
\begin{itemize}
\item[(1)] $\Gamma(t)$ is increasing in time, i.e., $\Gamma(t_1) \subseteq \Gamma(t_2)$
for all $0 \le t_1 \le t_2 \le 1$ (irreversibility of the process);
\item[{}]
\item[(2)] $\Es(g(t),\Gamma(t)) \le \Es(g(t),\Gamma)$ for all cracks $\Gamma$ such that
$\cup_{s<t} \Gamma(s) \subseteq \Gamma$ (equilibrium condition);
\item[{}]
\item[(3)] the total energy $\Es(g(t),\Gamma(t))$ is absolutely continuous in time
and (conservation of energy)
$$
\frac{d}{dt} \Es(g(t),\Gamma(t))= 2 \mu \int_\Om Eu(t) E\dot{g}(t)\,dx
+\lambda \int_\Om {\rm div}\,u(t) \,{\rm div}\,\dot{g}(t) \,dx.
$$
\end{itemize}
A precise mathematical formulation of this model has been given by
Dal Maso and Toader \cite{DMT} in the case of {\it anti-planar shear} in dimension two assuming
that the fractures are compact sets with a finite number of connected components.
Recently Francfort and Larsen \cite{FL}, using the framework of $SBV$-functions,
proved the existence of a quasi-static growth of brittle fractures in the case
of {\it generalized anti-planar shear} and without assumptions on the structure of the fractures
which are dealt with the set of jumps of the displacements.
To be precise, they consider as elastic body an infinite cylinder whose section $\Om \subseteq \R^N$
is subject to a displacement $u \in SBV(\Om)$ in the direction orthogonal to $\Om$.
The crack at time $t$ on the section $\Om$ is defined as
$$
\Gamma(t):=\bigcup_{s<t} \left[ S_{u(s)} \cup
(\partial_D \Om \cap \{u(s) \not= g(s)\})
\right],
$$
where $S_u$ denotes the set of jumps of $u$. Moreover the pair $(u(t), \Gamma(t))$
is such that:
\begin{itemize}
\item[(a)]
for all $v \in SBV(\Om)$
\begin{equation}
\label{intrfl1}
\int_\Om |\nabla u(t)|^2\,dx + \hs^{N-1}(\Gamma(t)) \le
\int_\Om |\nabla v|^2\,dx + \hs^{N-1}(S_v \cup
(\partial_D \Om \cap \{v \not= g(t)\})
\cup \Gamma(t));
\end{equation}
\item[{}]
\item[(b)]
the total energy $\Es(t):= \int_\Om |\nabla u(t)|^2\,dx+ \hs^{N-1}(\Gamma(t))$ is absolutely
continuous and
\begin{equation}
\label{intrfl2}
\Es(t)=\Es(0) + 2\int_0^t \int_\Om \nabla u(\tau) \nabla \dot{g}(\tau) \,dx \,d\tau.
\end{equation}
\end{itemize}
The aim of this paper is to discretize the model using a suitable finite
element method and to a give a rigorous proof of its convergence to a
quasi-static evolution in the sense of Francfort and Larsen.
We restrict our analysis to a two dimensional setting considering only a
polygonal reference configuration $\Om \subseteq \R^2$.
\par
The discretization of the domain $\Om$ is carried out, following \cite{N},
considering two parameters $\eps>0$ and $a \in ]0,\frac{1}{2}[$ .
We consider a regular triangulation $\rb_\eps$ of size $\eps$
of $\Om$, i.e. we assume that there exist two constants $c_1$ and $c_2$ so that
every triangle $T \in \rb_\eps$ contains a ball of radius $c_1 \eps$
and is contained in a ball of radius $c_2 \eps$. In order to treat the boundary data,
we assume also that $\partial_D \Om$ is composed of edges of $\rb_\eps$.
On each edge $[x,y]$ of $\rb_\eps$ we consider a point $z$ such that $z=tx+(1-t)y$ with
$t \in [a, 1-a]$. These points are called {\it adaptive vertices}. Connecting
together the adaptive vertices, we divide every $T \in \rb_\eps$ into four triangles.
We take the new triangulation $\tb$ obtained after this division as the discretization
of $\Om$. The family of all such triangulations is denoted by $\ts_{\eps,a}(\Om)$.
\par
The discretization of the energy functional is obtained restricting
the total energy to the family of functions $u$
which are affine on the triangles of some triangulation $\tb(u) \in \ts_{\eps,a}(\Om)$
and are allowed to jump across the edges of $\tb(u)$. We indicate this
space by ${\mathcal A}_{\eps,a}(\Om)$. The boundary data is assumed to belong
to the space $\afeom$ of continuous functions which are affine on every triangle $T \in
\rb_\eps$.
\par
Given the boundary data $g \in W^{1,1}([0,1],H^1(\Om))$ with $g(t) \in \afeom$ for all
$t \in [0,1]$,
we divide $[0,1]$ into subintervals $[t^\delta_i,t^\delta_{i+1}]$
of size $\delta>0$ for $i=0, \ldots, N_\delta$, we set
$g^\delta_i=g(t^\delta_i)$, and
for all $u \in \afeaom$ we indicate by $S_D^{g^\delta_i}(u)$ the edges
of the triangulation $\tb(u)$ contained in $\partial_D \Om$ on which $u \not= g^\delta_i$.
Using a variational argument we construct a {\it discrete evolution}
$\{u^{\delta,i}_{\eps,a}\,:\,i=0, \ldots,N_\delta\}$ such that
$u^{\delta,i}_{\eps,a} \in \afeaom$ for all $i=0, \ldots, N_\delta$, and
such that considering the {\it discrete fracture}
\begin{equation*}
\Gamma^{\delta,i}_{\eps,a}:=
\bigcup_{r=0}^i \big[ S_{u^{\delta,r}_{\eps,a}} \cup
S_D^{g^\delta_r}(u^{\delta,r}_{\eps,a}) \big],
\end{equation*}
the following {\it unilateral minimality property} holds:
\begin{equation}
\label{pieceminintr}
\int_\Om |\nabla u^{\delta,i}_{\eps,a}|^2\,dx \le \int_\Om |\nabla v|^2\,dx +
\hs^1\left( \big( S_v \cup S_D^{g^\delta_i}(v) \big) \setminus
\Gamma^{\delta,i-1}_{\eps,a}\right).
\end{equation}
Moreover we get suitable estimates for the discrete total energy
$$
\Es^{\delta,i}_{\eps,a}:=\|\nabla u^{\delta,i}_{\eps,a}\|^2_{L^2(\Om;\R^2)}+
\hs^1 \left( \Gamma^{\delta,i}_{\eps,a}\right).
$$
The definition of the discrete fracture ensures that
$\Gamma^{\delta,i}_{\eps,a} \subseteq \Gamma^{\delta,j}_{\eps,a}$ for all
$i \le j$, recovering in this discrete setting the irreversibilty of the growth
given in $(1)$.
The minimality property \eqref{pieceminintr} is the reformulation in the finite element
space of the equilibrium condition $(2)$.
\par
In order to perform the asymptotic analysis of the {\it discrete evolution}
$\{u^{\delta,i}_{\eps,a}\,:\,i=0, \ldots,N_\delta\}$, we make
the piecewise constant interpolation in time
$u^{\delta}_{\eps,a}(t)=u^{\delta,i}_{\eps,a}$ and
$\Gamma^{\delta}_{\eps,a}(t)=\Gamma^{\delta,i}_{\eps,a}$ for all $t^\delta_i \le t <t^\delta_{i+1}$.
The main result of the paper is the following theorem.

\begin{theorem}
\label{mainthm}
Let $g \in W^{1,1}([0,1], H^1(\Om))$ be such that $\|g(t)\|_\infty \le C$
for all $t \in [0,1]$ and let $g_\eps \in W^{1,1}([0,1], H^1(\Om))$ be
such that $\|g_\eps(t)\|_\infty \le C$, $g_\eps(t) \in \afeom$ for all $t \in [0,1]$
and
\begin{equation}
\label{bdryconv}
g_\eps \to g \quad \mbox { strongly in }W^{1,1}([0,1], H^1(\Om)).
\end{equation}
Given the discrete evolution $\{t \to u^{\delta}_{\varepsilon,a}(t)\}$
relative to the boundary data $g_\eps$, let $\Gamma^\delta_{\eps,a}$ and
$\Es^\delta_{\eps,a}$ be the associated fracture and total energy.
\par
Then there exist $\delta_n \to 0$, $\varepsilon_n \to 0$, $a_n \to 0$,
and a quasi-static evolution
$\{t \to (u(t),\Gamma(t)),\,t \in [0,1]\}$ relative to the boundary data $g$,
satisfying \eqref{intrfl1} and \eqref{intrfl2}, and such that setting
$u_n:=u^{\delta_n}_{\varepsilon_n, a_n}$,
$\Gamma_n:=\Gamma^{\delta_n}_{\varepsilon_n, a_n}$,
$\Es_n:=\Es^{\delta_n}_{\varepsilon_n, a_n}$, the following hold:
\begin{itemize}
\item[(a)] if $\ns$ is the set of discontinuities of $\hs^1(\Gamma(\cdot))$,
for all $t \in [0,1] \setminus \ns$ we have
\begin{equation}
\label{mainconvgradt}
\nabla u_n(t) \to \nabla u(t) \quad \mbox{ strongly in }L^2(\Om;\R^2)
\end{equation}
and
\begin{equation}
\label{mainconvjumpt}
\lim_n \hs^1(\Gamma_n(t))= \hs^1(\Gamma(t));
\end{equation}
\item[{}]
\item[(b)] for all $t \in [0,1]$ we have
\begin{equation}
\label{mainconvent}
\lim_n \Es_n(t)=\Es(t).
\end{equation}
\end{itemize}
\end{theorem}

We conclude that we have the convergence of the total energy at each time
$t \in [0,1]$, and the separate convergence of bulk and surface energy
for all $t \in [0,1]$ except a countable set.
\par
In order to prove Theorem \ref{mainthm}, we proceed in two steps.
Firstly, we fix $a$ and let
$\delta \to 0$ and $\eps \to 0$. We obtain an evolution
$\{t \to u_a(t)\,:\,t \in [0,1]\}$ such that
$\nabla u^{\delta}_{\eps,a}(t) \to \nabla u_a(t)$ strongly in $L^2(\Om;\R^2)$
for all $t$ up to a countable set and such that the following minimality
property holds: for all $v \in SBV(\Om)$
\begin{equation}
\label{aminimality}
\int_\Om |\nabla \vuat|^2\,dx \le \int_\Om |\nabla v|^2\,dx+ \mu(a)
\hs^1 \left( \big( S_v \cup (\partial_D \Om \cap \{v \not= g(t)\}) \big)
\setminus \Gamma_a(t) \right),
\end{equation}
where $\mu:]0,\frac{1}{2}[ \to ]0,+\infty[$ is a function independent of $\eps$ and
$\delta$, such that $\mu \ge 1$, $\lim_{a \to 0} \mu(a)=1$ and
$\Gamma_a(t):=\bigcup_{s \le t, s \in D} S_{\vuas} \cup (\partial_D \Om \cap \{u_a(s) \not= g(s)\})$.
The minimality property \eqref{aminimality} takes into account possible
anisotropies that could be generated as $\delta$ and $\eps \to 0$:
in fact, since $a$ is fixed, we have that the angles of the triangles in $\ts_{\eps,a}(\Om)$
are between fixed values (determined by $a$), and so fractures with
certain directions cannot be approximated in length.
In the second step, we let $a \to 0$ and determine from $\{t \to u_a(t)\,:\,t \in [0,1]\}$
a quasi-static evolution $\{t \to u(t)\,:\,t \in [0,1]\}$ in the sense of Francfort
and Larsen. Then, using a diagonal argument, we find sequences $\delta_n \to 0$,
$\eps_n \to 0$, and $a_n \to 0$ satisfying Theorem \ref{mainthm}.
\par
The main difficulties arise in the first part of our analysis, namely when
$\delta,\eps \to 0$. The convergence $u^\delta_{\eps,a}(t) \to u_a(t)$ in
$SBV(\Om)$ for $t \in D \subseteq [0,1]$ countable and dense is easily
obtained by means of Ambrosio's Compactness Theorem. The minimality property
\eqref{aminimality} derives from its discrete version \eqref{pieceminintr}
using a variant of Lemma 1.2 of \cite{FL}: given $v \in SBV(\Om)$, we
construct $v^\delta_{\eps,a} \in {\mathcal A}_{\eps,a}(\Om)$ such that
\begin{equation}
\label{intrconv1}
\nabla v^\delta_{\eps,a} \to \nabla v
\quad \mbox{ strongly in }L^2(\Om;\R^2)
\end{equation}
and
\begin{multline}
\label{intrconv2}
\limsup_{\delta,\eps \to 0}
\hs^1 \left[ \big( S_{v^\delta_{\eps,a}} \cup
S_D^{g_\eps^\delta(t)}(v^\delta_{\eps,a}) \big)
\setminus \Gamma^\delta_{\eps,a}(t) \right] \le \\
\le \mu(a)
\hs^1 \left[ \big( S_v \cup
\big( \partial_D \Om \cap \{v \not= g(t) \} \big) \big)
\setminus \Gamma_a(t) \right],
\end{multline}
where $g_\eps^\delta(t):=g_\eps(t^\delta_i)$ for $t^\delta_i \le t <t^\delta_{i+1}$.
The main difference with respect to Lemma 1.2 of \cite{FL} is that we have
to find the approximating functions $v^\delta_{\eps,a}$ in the finite
element space ${\mathcal A}_{\eps,a}(\Om)$. This can be regarded
as an interpolation problem, so we try to construct
triangulations $\tb_\eps \in \ts_{\eps,a}(\Om)$ adapted to $v$ in order to obtain
\eqref{intrconv1} and \eqref{intrconv2}.
In all the geometric operations involved, we need to avoid
degeneration of the triangles of $\tb(u^\delta_{\eps,a}(t))$ which is
guaranteed from the fact that $a$ is constant: this is the principal reason
to keep $a$ fixed in the first step.
A second difficulty arises when $u_a(\cdot)$ is extended from $D$ to the entire
interval $[0,1]$: indeed it is no longer clear whether
$\nabla u^\delta_{\eps,a}(t) \to \nabla u_a(t)$ for $t \not\in D$. Since
the space ${\mathcal A}_{\eps,a}(\Om)$ is not a vector space, we cannot
provide an estimate on
$\|\nabla u^\delta_{\eps,a}(t)-\nabla u^\delta_{\eps,a}(s)\|$ with $s \in D$
and $s<t$: we thus cannot expect to recover the convergence at time $t$ from
the convergence at time $s$. We overcome this difficulty observing that
$\nabla u^\delta_{\eps,a}(t) \to \nabla \tilde{u}_a$ with $\tilde{u}_a$ satisfying a minimality
property similar to \eqref{aminimality} and then proving
$\nabla \tilde{u}_a=\nabla u_a(t)$ by a uniqueness argument for the gradients of the
solutions.
\par
The plan of the paper is the following. In Section \ref{prel} we give the basic
definitions and prove some auxiliary results. In Section \ref{devol}, we
prove the existence of a discrete evolution. In Section
\ref{convres} we prove the convergence of the discrete evolution to a quasi-static
evolution of brittle fractures in the sense of Francfort and Larsen. The proof
of minimality property \eqref{aminimality} requires a careful analysis to which is
dedicated Section \ref{secmin}. In Section \ref{remark} we show that the arguments
of Section \ref{convres} can be used to improve the convergence results
for the discrete in time approximation considered in \cite{FL}.

\section{Preliminaries}
\label{prel}

In this section we state the notation and prove some preliminaries employed in
the rest of the paper.

\vskip10pt
{\it Basic notation.}
We will employ the following basic notation:
\begin{itemize}
\item[-] $\Om$ is a polygonal open subset of $\R^2$;
\item[-] $L^p(\Om;\R^m)$
with $1 \le p < +\infty$ and $m \ge 1$ is the usual Lebesgue space of
$p$-summable $\R^m$-valued functions, and $L^p(\Om):=L^p(\Om;\R)$;
\item[-] for all $k \ge 1$ and $1 \le p \le +\infty$,
$W^{k,p}(\Om)$ is the usual Sobolev space of functions
in $L^p(\Om)$ with distributional derivatives of order $1, \ldots, k$ in
$L^p(\Om)$; we will write $H^k(\Om)$ for $W^{k,2}(\Om)$;
\item[-] if $u \in W^{k,p}(\Om)$, $\nabla u$ is its gradient;
\item[-] $\hs^{1}$ is the one-dimensional Hausdorff measure;
\item[-] $\|\cdot\|_{\infty}$ denotes the sup-norm;
\item[-] if $f \in L^2(\Om;\R^m)$, $\|f\|$ denotes the $L^2$-norm
of $f$;
\item[-] for all $A \subseteq \R^2$, $|A|$ denotes the Lebesgue measure of $A$;
\item[-] if $\mu$ is a measure on $\R^2$ and $A$ is a Borel subset of $\R^2$,
$\mu \res A$ denotes the restriction of $\mu$ to $A$, i.e.
$(\mu \res A)(B):=\mu(B \cap A)$ for all Borel sets $B \subseteq \R^2$;
\item[-] if $\sigma \in ]0,+\infty[$, $o_\sigma$ is such that
$\lim_{\sigma \to 0^+} o_\sigma=0$.
\end{itemize}

\vskip20pt
{\it Special functions of bounded variation.}
For the general theory of functions of bounded variation, we refer to \cite{AFP}; here
we recall some basic definitions and theorems we need in the sequel.
Let $A$ be an open subset of $\R^N$, and let $u: A \to \R^n$.
We say that $u \in BV(A;\R^n)$ if $u \in L^1(A;\R^n)$, and its distributional
derivative is a vector-valued Radon measure on $A$.
We say that $u \in SBV(A;\R^n)$ if $u \in BV(A;\R^n)$ and its distributional
derivative can be represented as
$$
Du(A)= \int_A \nabla u(x) \,dx+ \int_{A \cap S_u} (u^+(x)-u^-(x)) \otimes \nu_x
\,d\hs^{N-1}(x),
$$
where $\nabla u$ denotes the approximate gradient of $u$, $S_u$ denotes the set of
approximate jumps of $u$, $u^+$ and $u^-$ are the traces of $u$ on $S_u$, $\nu_x$
is the normal to $S_u$ at $x$, and $\hs^{N-1}$ is the $(N-1)$-dimensional Hausdorff measure.
The space $SBV(A;\R^n)$ is called the space of
{\it special functions of bounded variation}. Note that if $u \in SBV(A;\R^n)$, then
the singular part of $Du$ is concentrated on $S_u$ which turns out to be countably
$\hs^{N-1}$-rectifiable.
\par
The space $SBV$ is very useful when dealing with variational problems involving
volume and surface energies because of the following compactness and lower
semicontinuity result due to L.Ambrosio (see \cite{A1}, \cite{A2}, \cite{A3}, \cite{AFP}).

\begin{theorem}
\label{SBVcompact}
Let $A$ be an open and bounded subset of $\R^N$, and let $(u_k)$ be a sequence in
$SBV(A;\R^n)$. Assume that there exists $q>1$ and $c \ge 0$ such that
$$
\int_A |\nabla u_k|^q \,dx+ \hs^{N-1}(S_{u_k})+ ||u_k||_\infty \le c
$$
for every $k \in \N$. Then there exists a subsequence $(u_{k_h})$ and a function
$u \in SBV(A;\R^n)$ such that
\begin{align}
\label{sbvconv}
\nonumber
u_{k_h} \to u \quad {strongly \; in}\; L^1(A;\R^n); \\
\nabla u_{k_h} \weak \nabla u \quad {weakly \; in}\; L^1(A;M^{N \times n}); \\
\nonumber
\hs^{N-1}(S_u) \le \liminf_h \hs^{N-1}(S_{u_{k_h}}).
\end{align}
\end{theorem}

In the rest of the paper, we will say that $u_k \to u$ in $SBV(A;\R^n)$ if
$u_k$ and $u$ satisfy \eqref{sbvconv}.
It will also be useful the following fact which can be derived from
Ambrosio's Theorem: if $u_k \to u$ in $SBV(A;\R^n)$ and if
$\hs^{N-1} \res S_{u_k} \weakst \mu$ weakly-star in the
sense of measures, then $\hs^{N-1} \res S_u \le \mu$ as measures.
We will set $SBV(A):=SBV(A;\R)$.

\vskip20pt
{\it Quasi-static evolution of brittle fracture.}
Let $\Om$ be an open bounded subset of $\R^N$ with Lipschitz boundary, and let
$\partial_D \Om$ be a subset of $\partial \Om$ open in the relative topology.
Let $g:[0,1] \to H^1(\Om)$ be absolutely continuous; we indicate
the gradient of $g$ at time $t$ by $\nabla g(t)$, and the time
derivative of $g$ at time $t$ by $\dot{g}(t)$. The main result of \cite{FL} is
the following theorem.

\begin{theorem}
\label{qse}
There exists a crack $\Gamma(t) \subseteq \Omb$ and a field $u(t) \in SBV(\Om)$ such
that:
\begin{itemize}
\item[(a)] $\Gamma(t)$ increases with $t$;
\item[{}]
\item[(b)] $u(0)$ minimizes
$$
\int_\Om |\nabla v|^2 \,dx +
\hs^{N-1}( S_v \cup
\{x \in \partial_D \Om: v(x) \not= g(0)(x)\})
$$
among all $v \in SBV(\Om)$ (inequalities on $\partial_D \Om$ are intended for the
traces of $v$ and $g$);
\item[{}]
\item[(c)]
for $t>0$, $u(t)$ minimizes
$$
\int_\Om |\nabla v|^2 \,dx
+\hs^{N-1} \left( \left[ S_v \cup
\{x \in \partial_D \Om: v(x) \not= g(t)(x)\} \right] \setminus \Gamma(t)
\right)
$$
among all $v \in SBV(\Om)$;
\item[{}]
\item[(d)]
$S_{u(t)} \cup \{x \in \partial_D \Om\,:\, u(t)(x) \not= g(t)(x) \}
\subseteq \Gamma(t)$, up to a set of $\hs^{N-1}$-measure $0$.
\end{itemize}
Furthermore, the total energy
$$
\Es(t):= \int_\Om |\nabla u(t)|^2 \,dx +\hs^{N-1}( \Gamma(t))
$$
is absolutely continuous and satisfies
$$
\Es(t)=\Es(0)+2 \int_0^t \int_\Om \nabla u(\tau) \nabla \dot{g}(\tau) \,dx \,d\tau
$$
for every $t \in [0,1]$. Finally, for any countable, dense set $I \subseteq [0,1]$,
the crack $\Gamma(t)$ and the displacement $u(t)$ can be chosen so that
$$
\Gamma(t)= \bigcup_{\tau \in I, \tau \le t}
\left( S_{u(\tau)} \cup
\{x \in \partial_D \Om\,:\, u(\tau)(x) \not= g(\tau)(x) \} \right).
$$
\end{theorem}

The main tool in the proof of Theorem \ref{qse} is the following result
\cite[Theorem 2.1]{FL}, which is useful also in our analysis.

\begin{theorem}
\label{jumptransfer}
Let $\Omb \subseteq \Om'$, with $\partial \Om$ Lipschitz, and let for $r=1,\dots,i$
$(u^r_n)$ be a sequence in $SBV(\Om')$ such that
\begin{itemize}
\item[(1)] $S_{u^r_n} \subseteq \Omb$;
\item[{}]
\item[(2)] $|\nabla u^r_n|$ weakly converges in $L^1(\Om')$; and
\item[{}]
\item[(3)] $u^r_n \to u^r$ strongly in $L^1(\Om')$,
\end{itemize}
where $u^r \in BV(\Om')$ with $\hs^{N-1}(S_{u^r}) <\infty$. Then for every $\phi \in
SBV(\Om')$ with $\hs^{N-1}(S_\phi) <\infty$ and $\nabla \phi \in L^q(\Om;\R^N)$ for some $q \in [1,+\infty[$, 
there exists a sequence $(\phi_n)$ in
$SBV(\Om')$ with $\phi_n \equiv \phi$ on $\Om' \setminus \Omb$ such that
\begin{itemize}
\item[(a)] $\phi_n \to \phi$ strongly in $L^1(\Om')$;
\item[{}]
\item[(b)] $\nabla \phi_n \to \nabla \phi$ strongly in $L^q(\Om')$; and
\item[{}]
\item[(c)]
$\hs^{N-1} \left( [S_{\phi_n} \setminus \bigcup_{r=1}^i S_{u^r_n}] \setminus
[S_{\phi} \setminus \bigcup_{r=1}^i S_{u^r}] \right) \to 0$.
\end{itemize}
In particular
\begin{equation}
\label{jt}
\limsup_n \hs^{N-1} \left( S_{\phi_n} \setminus \bigcup_{r=1}^i S_{u^r_n} \right) \le
\hs^{N-1} \left( S_{\phi} \setminus \bigcup_{r=1}^i S_{u^r} \right).
\end{equation}
\end{theorem}

\vskip 20pt
{\it Hausdorff metric on compact sets.}
Let $A \subseteq \R^2$ be open and bounded, and let $\ks(\overline{A})$
be the set of all compact subsets of $\overline{A}$.
$\ks(\overline{A})$ can be endowed by the
Hausdorff metric $d_H$ defined by
$$
d_H(K_1,K_2) := \max \left\{ \sup_{x \in K_1} {\rm dist}(x,K_2), \sup_{y \in
K_2} {\rm dist}(y,K_1)\right\},
$$
with the conventions ${\rm dist}(x, \emptyset)= {\rm diam}(A)$ and $\sup
\emptyset=0$, so that $d_H(\emptyset, K)=0$ if $K=\emptyset$ and
$d_H(\emptyset,K)={\rm diam}(A)$ if $K \not=\emptyset$. It turns out that
$\ks(\overline{A})$ endowed with the Hausdorff metric is a compact space
(see e.g. \cite{Ro}).

\vskip20pt
{\it Triangulations.}
Let $\Om \subseteq \R^2$ be a polygonal set and let us fix two positive
constants $0<c_1<c_2$. By a {\it regular triangulation} of $\Om$ of
size $\eps$ we intend a finite family of (closed) triangles $T_i$ such that $\Omb=\bigcup_i T_i$,
$T_i \cap T_j$ is either empty or equal to a common edge or to a common vertex,
and each $T_i$ contains a ball of diameter $c_1 \eps$ and is contained in a ball of
diameter $c_2 \eps$.
\par
We indicate by $\rs_\eps(\Om)$ the family of all regular triangulations of
$\Om$ of size $\eps$.
It turns out that there exist $0<\vartheta_1 < \vartheta_2 <\pi$
such that for all $T$ belonging to a triangulation
$\tb \in \rs_\eps(\Om)$, the inner angles of $T$ are between
$\vartheta_1$ and $\vartheta_2$. Moreover, every edge of $T$ has length greater than
$c_1 \eps$ and lower than $c_2 \eps$.
\par
Let us fix a triangulation $\rb_\eps \in \rs_\eps(\Om)$ for all $\eps>0$
and let $a \in ]0,\frac{1}{2}[$.
Let us consider a new triangulation $\tb$ nested in $\rb_\eps$ obtained
dividing each $T \in \rb_\eps$ into four triangles taking over every edge $[x,y]$ of $T$
a knot $z$ which satisfies
$$
z=tx+(1-t)y, \quad \quad t \in [a,1-a].
$$
We will call these new vertices {\it adaptive}, the triangles obtained joining these points
{\it adaptive triangles}, and their edges {\it adaptive edges} (see Fig.1).
\begin{center}
\psfig{figure=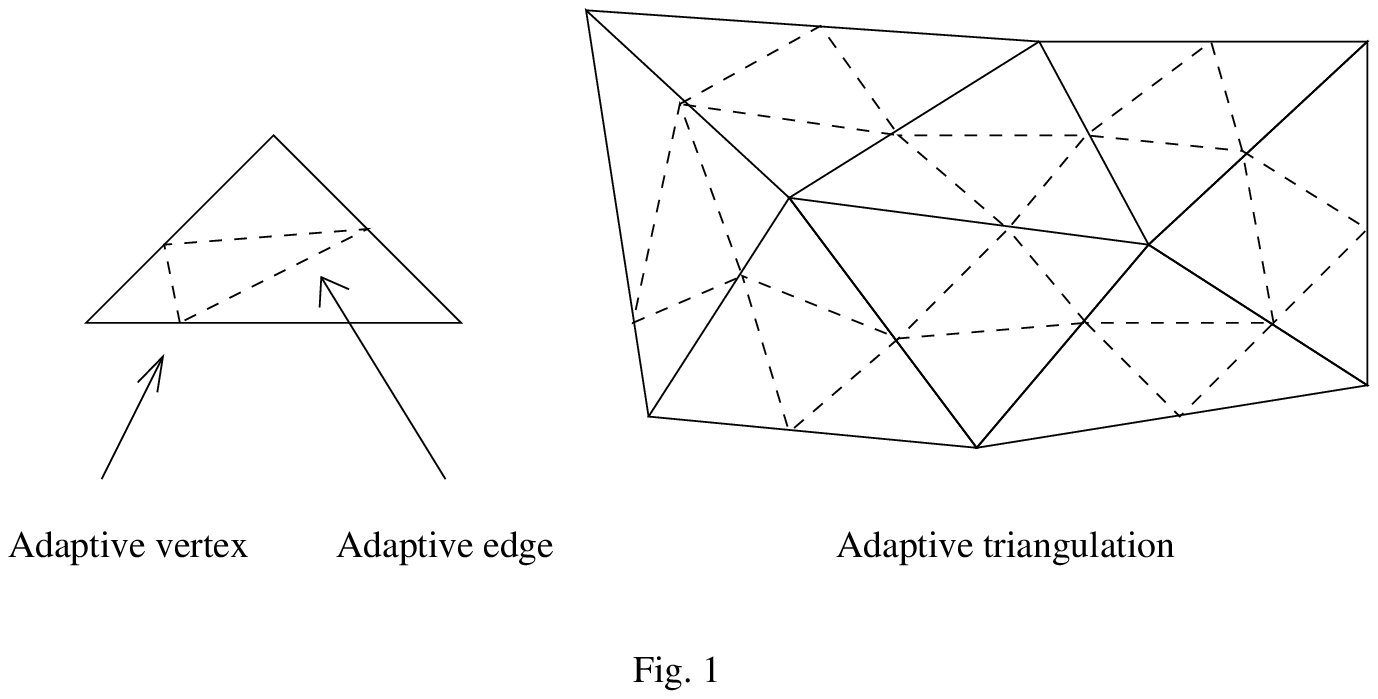}
\end{center}
We denote by $\ts_{\eps,a}(\Om)$ the set of all triangulations $\tb$ constructed in this way.
Note that for all $\tb \in \ts_{\eps,a}(\Om)$ there exists $0<c_1^a<c_2^a<+\infty$ such that
every $T_i \in \tb$ contains a ball of diameter $c_1^a \eps$ and is contained in a ball
of diameter $c_2^a \eps$. Then there exist $0<\vartheta_1^a < \vartheta_2^a <\pi$
such that for all triangles $T$ belonging to a triangulation 
$\tb \in \ts_{\eps,a}(\Om)$, the inner angles of $T$ are between
$\vartheta_1^a$ and $\vartheta_2^a$. Moreover, every edge of $T$ has length greater than
$c_1^a \eps$ and lower than $c_2^a \eps$.
\par
We will often use the following {\it interpolation estimate} (see \cite[Theorem 3.1.5]{Cia}).
If $u \in W^{2,2}(\Om)$ and $T \in \rb_\eps$, let $u_T$ denote the affine interpolation
of $u$ on $T$. We have that there exists $\cint$ depending only on $c_1,c_2$ such that
\begin{equation}
\label{interp1}
\|u_T -u\|_{W^{1,2}(T)} \le \cint \eps \|u\|_{W^{2,2}(T)}.
\end{equation}
Estimate \eqref{interp1} holds also for
$\tb \in \ts_{\eps,a}(\Om)$: in this case $K$ depends on $a$.

\vskip20pt
{\it Some elementary constructions.}
The following lemmas will be used in Section 4.

\begin{lemma}
\label{simplecurve}
Let $\tb \in \ts_{\eps,a}(\Om)$, and let $l \subseteq \Om$
be a segment with extremes $p,\,q$ belonging to edges of $\tb$.
There exists a polyhedral curve $\Gamma$ with extremal points $p$ and $q$
(see Fig.2) such that $\Gamma$ is contained in the union of the edges of those
$T \in \tb$ with $T \cap l \neq \emptyset$, and such that the following
properties  hold:
\begin{itemize}
\item[(1)]
$\Gamma = \gamma_p \cup \gamma \cup \gamma_q$, where $\gamma$ is union of edges
of $\tb$ and $\gamma_p$, $\gamma_q$ are segments containing $p$ and $q$ respectively,
and each one is contained in an edge of $\tb$;
\item[]
\item[(2)]
there exists a constant $c$ independent of $\epsilon$ (but depending on $a$) such that
$$
\hs^1(\Gamma)\le c \, \hs^1(l).
$$
\end{itemize}
\end{lemma}

\begin{center}
\psfig{figure=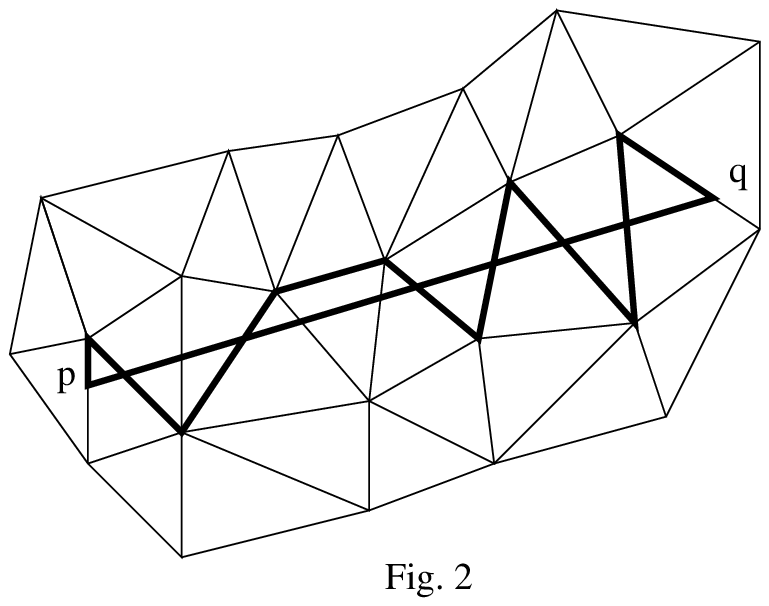}
\end{center}

\begin{proof}
Let $\{T_1, \ldots T_k\}$ be the family of
triangles in $\tb$ such that the intersection with $l$ is a segment with positive
length. For every integer $1\le i \le k$,
let $l_i:=T_i\cap l.$  If $l_i$ is an edge of $T_i$, we set $D_i= T_i$. Otherwise
let $D_i$ be a connected component of $T_i \setminus l_i$ such that
$|D_i| \leq \frac{1}{2} |T_i|.$
We claim that there exists a constant $c>0$ independent of $\epsilon$ such that
\begin{equation}
\label{inelung}
\hs^1 (\partial D_i) \leq c\, \hs^1(l_i).
\end{equation}
We have to analyze two possibilities, namely $D_i$ is a triangle, or $D_i$ is a
trapezoid. Suppose that $D_i$ is a triangle and that $m_i$ is an edge of $D_i$.
Let $\alpha$ be the angle of $D_i$ opposite to
$l_i$. It is easy to prove that $\hs^1(l_i) \ge \hs^1(m_i) \sin\alpha$,
and so
$$
\hs^1(l_i) \ge \frac{1}{3} \sin\alpha \hs^1(\partial D_i).
$$
Since $\vartheta_1^a\le\alpha\le\vartheta_2^a$,  $\sin\alpha$ is uniformly bounded
from below, and hence inequality (\ref{inelung}) follows.
If $D_i$ is a trapezoid, since $|D_i| \leq \frac{1}{2} |T_i|$,
it follows that $T_i\setminus D_i$ is a triangle such that
its edges different from $l_i$ have length greater than
$\frac{1}{2} c_1^a \eps$. Let $\alpha$ be the inner angle of $T_i\setminus D_i$
opposite to $l_i$. We have that
$$
\hs^1(l_i) \ge \frac{1}{2}\sin\alpha c_1^a {\eps} \ge
\frac{1}{2} \sin\alpha \frac {c_1^a}{c_2^a}\frac{1}{4} \hs^1(\partial D_i).
$$
Since $\vartheta_1^a\le\alpha\le\vartheta_2^a$, inequality (\ref{inelung}) follows.
\par
By \eqref{inelung}, we deduce that
\begin{equation*}
\hs^1 \left( \bigcup_{i=1}^k\partial D_i \right) \leq c\, \hs^1(l);
\end{equation*}
moreover, since $\bigcup_{i=1}^k\partial D_i$ is arcwise connected and contains
$p,q$, we conclude that there exists a curve $\Gamma \subseteq
\bigcup_{i=1}^k\partial D_i$ which satisfies the thesis.
\end{proof}

\begin{lemma}
\label{guscio}
There exists a constant $c>0$ such that for
every segment $l \subseteq \Om$ there exists $\eps_0$ with the following property:
for every $\varepsilon \leq \varepsilon_0$, setting
$\rs(l):=\{T\in \rb_\eps:\, T \cap l \neq \emptyset\}$, we have
$$
\hs^1(\partial \rs(l))\leq c \hs^1(l).
$$
\end{lemma}

\begin{proof}
Let $ \ns_{\eps}(l):=\{ x \in \Om:\, {\rm dist}(x,l) \leq c_2 \eps\}$.
We have that
$|\ns_{\eps}(l)| = \hs^1(l)c_2\eps+ \pi c_2^2\eps^2$,
and hence there exists a positive constant $\eps_0$ such that,
for every $\eps \leq \eps_0$, we have that
$$
|\ns_{\eps}(l)| \leq 2 \hs^1(l) c_2 \eps.
$$
We have that $\rs(l) \subseteq \ns_{\eps}(l),$ and
$$
\sharp \rs(l) \le \frac{4}{c_1^2\pi^2} \frac{|\ns_{\eps}(l)|}{\eps^2},
$$
where $\sharp \rs(l)$ denotes the number of triangles of $\rs(l)$.
Then, we have
$$
\hs^1(\partial \rs(l)) \leq 3 c_2 \eps\, \sharp\rs(l) \leq
3 c_2 \eps   \frac{4}{c_1^2\pi^2}
\frac{|\ns_{\eps}(l)|}{\eps^2} \leq
3 c_2^2 \frac{4}{c_1^2\pi^2} 2\hs^1(l),
$$
and so the proof is concluded.
\end{proof}

\vskip20pt
{\it A density result.}
Let $A \subseteq \R^2$ be open.
We say that $K \subseteq A$ is polygonal (with respect to $A$),
if it is the intersection of $A$ with the union of a finite number of closed segments.
The following density result is proved in \cite{C}.

\begin{theorem}
\label{piecedensity}
Assume that $\partial A$ is locally Lipschitz, and let $u \in SBV(A)$ such that
$u \in L^2(A)$, $\nabla u \in L^2(\Om;\R^2)$, and $\hs^{1}(S_u)<+\infty$. For
every $\eps>0$, there exists a function $v \in SBV(A)$ such that
\begin{itemize}
\item[(a)]
$S_v$ is essentially closed, i.e., $\hs^{1}( \overline{S_v} \setminus S_v)=0$;
\item[{}]
\item[(b)]
$\overline{S_v}$ is a polyhedral set;
\item[{}]
\item[(c)]
$v \in W^{k, \infty}(A \setminus \overline{S_v})$ for every $k \in \N$;
\item[{}]
\item[(d)]
$||v-u||_{L^2(A)} < \eps$;
\item[{}]
\item[(e)]
$||\nabla v- \nabla u||_{L^2(A; \R^2)} < \eps$;
\item[{}]
\item[(f)]
$|\hs^{1}(S_v)-\hs^{1}(S_u)| <\eps$.
\end{itemize}
\end{theorem}

Let $\partial_D \Om$ be a relatively open subset of $\partial \Om$ composed of edges lying in 
$\partial \Om$. Let us consider $\Om_D$ polygonal open bounded subset of $\R^2$ such that 
$\Om_D \cap \Om=\emptyset$ and $\partial \Om \cap \partial \Om_D= \partial_D \Om$ up to a finite 
number of vertices. We set $\Om':= \Om \cup \Om_D \cup \partial_D \Om$.
In Section 4, we will use the following result.

\begin{proposition}
\label{regularization}
Given $u \in SBV(\Om')$ with $u=0$ on $\Om' \setminus \Omb$ and $\hs^{N-1}(S_u)<+\infty$,
there exists $u_h \in SBV(\Om')$ such that
\begin{itemize}
\item[(a)] $u_h=0$ in $\Om' \setminus \Omb$;
\item[{}]
\item[(b)] $S_{u_h}$ is polyhedral,
$\overline{S_{u_h}} \subseteq \Om$ and $u_h \in 
W^{k,\infty}(\Om' \setminus \overline{S_{u_h}})$ for all $k$;
\item[{}]
\item[(c)] $u_h \to u$ strongly in $L^2(\Om')$ and
$\nabla u_h \to \nabla u$ strongly in $L^2(\Om';\R^2)$;
\item[{}]
\item[(d)] for all $A$ open subset of $\Om'$ with $\hs^1(\partial A \cap S_u)=0$,
we have
$$
\lim_h \hs^1(A \cap S_{u_h})=\hs^1(A \cap S_u).
$$
\end{itemize}
\end{proposition}

\begin{proof}
Using a partition of unity, we may prove the result in the case $\Om:=]-1,1[ \times ]0,1[$,
$\Om':=]-1,1[ \times ]-1,1[$, and $\partial_D \Om:=]-1,1[ \times \{0\}$.
We set $w_h(x,y):=u(x,y-h)$, and let $\varphi_h$ be a cut off function with $\varphi_h=1$ on
$]-1,1[ \times ]-1,\frac{h}{3}[$, $\varphi_h=0$ on $]-1,1[ \times ]\frac{h}{2},1[$, and
$||\nabla \varphi_h||_\infty \le \frac{7}{h}$. Let us set $v_h:=(1-\varphi_h)w_h$.
We have that $v_h=0$ in $\Om' \setminus \Omb$; moreover we have
$$
\nabla v_h=(1-\varphi_h) \nabla w_h-\nabla \varphi_h w_h.
$$
Since $\nabla \varphi_h w_h \to 0$ strongly in $L^2(\Om';\R^2)$, we have
$\nabla v_h \to \nabla u$ strongly in $L^2(\Om';\R^N)$. Finally, for all
$A$ open subset of $\Om'$ with $\hs^1(\partial A \cap S_u)=0$,
we have
$$
\lim_h \hs^1(A \cap S_{v_h})=\hs^1(A \cap S_u).
$$
In order to conclude the proof, let us apply Theorem \ref{piecedensity}
obtaining $\tilde{v}_h$ with polyhedral jumps in $\Om$ such that $\tilde{v}_h
\in W^{k,\infty}(\Om' \setminus \overline{S_{\tilde{v}_h}})$,
$||w_h-\tilde{v}_h||_{L^2(\Om)}+
||\nabla w_h -\nabla \tilde{v}_h||_{L^2(\Om;\R^2)} \le h^2$ and
$|\hs^{N-1}(S_{w_h}) -\hs^{N-1}(S_{\tilde{v}_h})| \le h$. If we set
$u_h:= \varphi_h g+(1-\varphi_h) \tilde{v}_h$, we obtain the thesis.
\end{proof}

\section{The discontinuous finite element approximation}
\label{devol}

In this section we construct a discrete approximation of quasi-static evolution
of brittle fractures in linearly elastic bodies: the discretization is done both in
space and time.
\par
From now on we suppose that $\Om$ is a polygonal open bounded subset
of $\R^2$, and that $\partial_D \Om \subseteq \partial \Om$ is open in the relative
topology.
For all $\varepsilon>0$, we fix a triangulation $\rb_\varepsilon \in \rs_\varepsilon(\Om)$,
and suppose that $\partial_D \Om$ is composed of edges of $\rb_\eps$ for
all $\eps$; we indicate the family of these edges by $\seub_\eps$.
\par
We consider the following discontinuous finite element space.
We indicate by $\afeaom$ the set of all $u$ such that there exists a triangulation
$\tb(u) \in \treaom$ nested in $\rb_\varepsilon$ with $u$ affine on every
$T \in \tb(u)$. For every $u \in \afeaom$, we write $\|\nabla u\|$
for the $L^2$-norm of $\nabla u$ and we indicate by $S_u$ the family of edges of
$\tb(u)$ inside $\Om$ across which $u$ is discontinuous. Notice that $u \in SBV(\Om)$
and that the notation is consistent with the usual one employed in the theory of
functions with bounded variation. Let us also denote by $\afeom$ the set of affine
functions in $\Om$ with respect to the triangulation $\rb_\eps$.
Finally, given any $g \in \afeom$, for all $u \in \afeaom$ we set
\begin{equation}
\label{jump*}
S_D^g(u):= \{\zeta \in \seub_\eps\,:\, u \not=g \mbox{ on } \zeta\},
\end{equation}
that is $S_D^g(u)$ denotes the edges at which the boundary condition is not
satisfied. Moreover we set
\begin{equation}
\label{jumpbis}
\Sg{g}{u}:= S_u \cup S_D^g(u)
\end{equation}
Let now consider $g \in W^{1,1}([0,1];H^1(\Om))$ with $g(t) \in \afeom$ for all
$t \in [0,1]$. Let $\delta>0$ and let $N_\delta$ be the largest integer such that
$\delta (N_\delta-1) < 1$; for $0\le i \le N_\delta-1$ we set $t_i^\delta:=i\delta$,
$t^\delta_{N_\delta}:=1$ and $g_i^\delta:=g(t_i^\delta)$.
The following proposition holds.

\begin{proposition}
\label{discrevol}
Let $\eps>0$, $a \in ]0,\frac{1}{2}]$ and $\delta>0$ be fixed. Then for all
$i=0, \ldots, N_\delta$ there exists $u^{\delta,i}_{\eps,a} \in \afeaom$
such that, setting
\begin{equation*}
\Gamma^{\delta,i}_{\eps,a}:=
\bigcup_{r=0}^i \Sg{g^\delta_r}{u^{\delta,r}_{\eps,a}},
\end{equation*}
the following hold:
\begin{itemize}
\item[]
\item[(a)] $\|u^{\delta,i}_{\eps,a}\|_\infty \le \|g^\delta_i\|_\infty$;
\item[{}]
\item[(b)] for all $v \in \afeaom$ we have
\begin{equation}
\label{MSdiscr}
\|\nabla u^{\delta,0}_{\eps,a}\|^2+
\hs^1 \left(
\Sg{g_0^\delta}{u^{\delta,0}_{\eps,a}} \right)
\le \|\nabla v\|^2 +
\hs^1\left( \Sg{g_0^\delta}{v} \right),
\end{equation}
and
\begin{equation}
\label{piecemin}
\|\nabla u^{\delta,i}_{\eps,a}\|^2 \le \|\nabla v\|^2 +
\hs^1\left( \Sg{g^\delta_i}{v} \setminus
\Gamma^{\delta,i-1}_{\eps,a}\right).
\end{equation}
\end{itemize}
\end{proposition}

\begin{proof}
The proof is carried out through a variational argument.
Let $u^{\delta,0}_{\eps,a}$ be a minimum of the following problem
\begin{equation}
\label{step0}
\min
\left\{ \|\nabla u\|^2+ \hs^1(\Sg{g^\delta_0}{u}) \right\}.
\end{equation}
We set $\Gamma^{\delta,0}_{\eps,a}:= \Sg{g^\delta_0}{u^{\delta,0}_{\eps,a}}$.
Recursively, supposing to have constructed $u^{\delta,i-1}_{\eps,a}$ and
$\Gamma^{\delta,i-1}_{\eps,a}$, let $u^{\delta,i}_{\eps,a}$ be a minimum for
\begin{equation}
\label{stepj}
\min \left\{ \|\nabla u\|^2 +
\hs^1 \left( \Sg{g^\delta_i}{u} \setminus
\Gamma^{\delta,i-1}_{\eps,a} \right) \right\}.
\end{equation}
We set $\Gamma^{\delta,i}_{\eps,a}:=
\Sg{g^\delta_i}{u^{\delta,i}_{\eps,a}} \cup \Gamma^{\delta,i-1}_{\eps,a}$.
We claim that problems \eqref{step0} and \eqref{stepj} admit a solution
$u^{\delta,i}_{\eps,a}$ such that
$\|u^{\delta,i}_{\eps,a}\|_{\infty} \le \|g^\delta_i\|_\infty$ for all
$i=0, \ldots, N_\delta$.
We prove the claim for problem \eqref{stepj}, the other case being similar.
Let $(u_n)$ be a minimizing sequence for problem \eqref{stepj}: since
$g^\delta_i$ is an admissible test function, we deduce that
for $n$ large
\begin{equation*}
\|\nabla u_n\|^2+ \hs^1 \left( \Sg{g^\delta_i}{u_n}
\setminus \Gamma^{\delta,i-1}_{\eps,a} \right) \le
\|\nabla g^\delta_i\|^2.
\end{equation*}
Moreover, we may modify $u_n$ in the following way.
If $\pi$ denotes the projection in
$\R$ over the interval $I:=[-\|g^\delta_i\|_\infty, \|g^\delta_i\|_\infty]$, let
$\tilde{u}_n \in \afeaom$ be defined on each $T \in \tb(u_n)$ as the affine
interpolation of the values $(\pi(u_n(x_1)),\; \pi(u_n(x_2)),\; \pi(u_n(x_3))$, where
$x_1, \;x_2$ and $x_3$ are the vertices of $T$. Note that by construction we have
for all $n$
\begin{equation*}
\|\tilde{u}_n\|_\infty \le \|g^\delta_i\|_\infty, \quad \quad
\|\nabla \tilde{u}_n\| \le \|\nabla u_n\|, \quad \quad
\Sg{g^\delta_i}{\tilde{u}_n} \subseteq
\Sg{g^\delta_i}{u_n},
\end{equation*}
so that $(\tilde{u}_n)$ is a minimizing sequence for problem
\eqref{stepj}. We conclude that it is not restrictive to assume
$\|u_n\|_\infty \le \|g^\delta_i\|_\infty$.
\par
Since $\tb(u_n) \in \treaom$, we have that the number of elements of
$\tb(u_n)$ is uniformly bounded. Up to a subsequence,
we may suppose that there exists an integer $k$ such that $\tb(u_n)$ has exactly $k$
elements $T_n^1, \ldots, T_n^k$. Using a diagonal argument we may suppose that, up
to a further subsequence, there exists $\tb=\{T^1, \ldots, T^k\} \in \treaom$ such that
$T_n^i \to T^i$ in the Hausdorff metric for all $i=1, \ldots, k$. Let us consider
$T^i \in \tb$, and let $\tilde{T}^i$ be contained in the interior of $T^i$. For $n$
large enough, $\tilde{T}^i$ is contained in the interior of $T^i_n$ and
$(u_n)_{|\tilde{T}^i}$ is affine with
$\int_{\tilde{T}^i} |\nabla u_n|^2 \,dx \le C$ with $C$ independent of $n$.
We deduce that there exists a function $u^i$ affine on $\tilde{T}^i$ such that up to a
subsequence $u_n \to u$ uniformly on $\tilde{T}^i$. Since $\tilde{T}^i$ is arbitrary,
it turns out that $u^i$ is actually defined on $T^i$ and
\begin{equation*}
\int_{T^i} |\nabla u^i|^2 \,dx \le \liminf_n \int_{T_n^i} |\nabla u_n|^2 \,dx.
\end{equation*}
Let $u \in \afeaom$ such that $u=u^i$ on $T^i$ for every $i=1,\dots,k$: we have
\begin{equation*}
\|\nabla u\|^2 \le \liminf_n \|\nabla u_n\|^2.
\end{equation*}
On the other hand, it is easy to see that
$\Sg{g^\delta_i}{u}$ is contained in the Hausdorff limit
of $\Sg{g^\delta_i}{u_n}$, and that
\begin{equation*}
\hs^1 \left(
\Sg{g^\delta_i}{u} \setminus \Gamma^{\delta,i-1}_{\eps,a}
\right)
\le
\liminf_n \hs^1 \left(
\Sg{g^\delta_i}{u_n}
\setminus \Gamma^{\delta,i-1}_{\eps,a}
\right).
\end{equation*}
We conclude that $u$ is a minimum point for the problem \eqref{stepj} with
$\|u\|_\infty \le \|g^\delta_i\|_\infty$. We have that point $(a)$ is proved.
\par
Concerning point $(b)$, by construction we get
\eqref{MSdiscr}; for $i \ge 1$ we have
\begin{equation*}
\|\nabla u^{\delta,i}_{\eps,a}\|^2 +
\hs^1 \left(
\Sg{g^\delta_i}{u^{\delta,i}_{\eps,a}}
\setminus \Gamma^{\delta,i-1}_{\eps,a} \right) \le
\|\nabla v\|^2 +
\hs^1 \left( \Sg{g^\delta_i}{v}
\setminus \Gamma^{\delta,i-1}_{\eps,a} \right)
\end{equation*}
for all $v \in \afeaom$, so that
\begin{equation*}
\|\nabla u^{\delta,i}_{\eps,a}\|^2 \le
\|\nabla v\|^2 +
\hs^1 \left( \Sg{g^\delta_i}{v} \setminus
\Gamma^{\delta,i-1}_{\eps,a} \right),
\end{equation*}
and this proves point $(b)$.
\end{proof}

\begin{remark}
\label{oss}
{\rm
For technical reasons due to the asymptotic analysis of the
discrete evolution $u^{\delta,i}_{\eps,a}$ when $\delta \to 0$,
$\eps \to 0$ and $a \to 0$, we define $u^{\delta,i}_{\eps,a}$
from $u^{\delta,i-1}_{\eps,a}$ through problem \eqref{stepj}
without requiring that the adaptive vertices determining
$\Gamma^{\delta,i-1}_{\eps,a}$ remain fixed. We just
penalize their possible changes if they are used to create new
fracture: in fact in this case, the surface energy increases at each change
of a quantity at least of order $a\eps$.
As a consequence, during the step by step minimization, it could
happen that some triangles $T \in \ts_{\eps,a}(\Om)$ contain the
fracture $\Gamma^{\delta,i}_{\eps,a}$ in their interior.
This is in contrast with the interpretation of the triangles as elementary
blocks for the elasticity problem, but being this situation penalized in the
minimization process, we expect that it occurs rarely.
}
\end{remark}

The following estimate is essential for the study of asymptotic behavior of
the discrete evolution.

\begin{proposition}
\label{dener}
If $(u^{\delta,i}_{\eps,a},\Gamma^{\delta,i}_{\eps,a})$ for $i=0,\dots,N_\delta$ satisfies condition
$(b)$ of Proposition \ref{discrevol}, setting $\Es^{\delta,i}_{\eps,a}:=\|\nabla u^{\delta,i}_{\eps,a}\|^2+
\hs^1 \left( \Gamma^{\delta,i}_{\eps,a}\right)$, we have
for $0 \le j \le i \le N_{\delta}$
\begin{equation}
\label{discrenergy}
\Es^{\delta,i}_{\eps,a} \le \Es^{\delta,j}_{\eps,a}+
2 \sum_{r=j}^{i-1} \int_{t^\delta_{r}}^{t^\delta_{r+1}}
\int_{\Om} \nabla u^{\delta,r}_{\eps,a} \nabla\dot{g}(\tau)
\,dx\,d\tau +o^\delta,
\end{equation}
where
\begin{equation}
\label{odelta}
o^\delta:= \left[
\max_{r=0,\ldots,N_\delta-1}
\int_{t^\delta_{r}}^{t^\delta_{r+1}} \|\dot{g}(\tau)\|_{H^1(\Om)} \,d\tau
\right]
\int_0^1 \|\dot{g}(\tau)\|_{H^1(\Om)}\,d\tau.
\end{equation}
\end{proposition}

\begin{proof}
For all $0 \le j \le N_\delta-1$,
by construction of $u^{\delta,j+1}_{\eps,a}$ we have that
\begin{multline*}
\|\nabla u^{\delta,j+1}_{\eps,a}\|^2 +
\hs^1 \left(
\Sg{g^\delta_{j+1}}{u^{\delta,j+1}_{\eps,a}}
\setminus \Gamma^{\delta,j}_{\eps,a} \right)
\le
\|\nabla u^{\delta,j}_{\eps,a}+ \nabla (g^\delta_{j+1} - g^\delta_{j})\|^2 =\\
=
\|\nabla u^{\delta,j}_{\eps,a}\|^2 +
2 \int_{\Om}
\nabla u^{\delta,j}_{\eps,a} \nabla (g^\delta_{j+1} - g^\delta_{j})
\,dx+||\nabla (g^\delta_{j+1}-g^\delta_{j})||^2.
\end{multline*}
Notice that
\begin{equation*}
\nabla (g^\delta_{j+1}-g^\delta_{j}) =
\int_{t^\delta_{j}}^{t^\delta_{j+1}}
\nabla \dot{g}(\tau) \,d\tau,
\end{equation*}
so that
\begin{multline}
\label{ineqi}
\|\nabla u^{\delta,j+1}_{\eps,a}\|^2 +
\hs^1 \left(
\Sg{g^\delta_{j+1}}{u^{\delta,j+1}_{\eps,a}}
\setminus \Gamma^{\delta,j}_{\eps,a} \right) \le \\
\le \|\nabla u^{\delta,j}_{\eps,a}\|^2 +
2 \int_{t^\delta_{j}}^{t_{j+1}^\delta} \int_{\Om}
\nabla u^{\delta,j}_{\eps,a} \nabla \dot{g}(\tau)\,dx \,d\tau
+e(\delta)\int_{t^\delta_{j}}^{t^\delta_{j+1}} \|\dot{g}(\tau)\|_{H^1(\Om)}
\,d\tau,
\end{multline}
where
\begin{equation*}
e(\delta):= \max_{r=0,\ldots,N_\delta-1}
\int_{t^\delta_{r}}^{t^\delta_{r+1}} \|\dot{g}(\tau)\|_{H^1(\Om)} \,d\tau.
\end{equation*}
From \eqref{ineqi}, we obtain that for all $0 \le j \le i \le N_\delta$
\begin{multline*}
||\nabla u^{\delta,i}_{\eps,a}||^2+ \hs^1(\Gamma^{\delta,i}_{\eps,a}) \le
||\nabla u^{\delta,j}_{\eps,a}||^2 + \hs^1(\Gamma^{\delta,j}_{\eps,a})+ \\
+2 \sum_{r=j}^{i-1} \int_{t^\delta_{r}}^{t^\delta_{r+1}} \int_{\Om}
\nabla u^{\delta,r}_{\eps,a} \nabla \dot{g}(\tau) \,dx\,d\tau
+e(\delta) \int_{t^\delta_{j}}^{t^\delta_{i}}
||\dot{g}(\tau)||_{H^1(\Om)} \,d\tau,
\end{multline*}
and so the proof of point $(c)$ is complete choosing
$$
o^\delta:=e(\delta) \int_0^1 \|\dot{g}(\tau)\|_{H^1(\Om)}\,d\tau.
$$
\end{proof}

\section{The convergence result}
\label{convres}
This section is devoted to the proof of Theorem \ref{mainthm}.
As in Section \ref{devol}, let $\Om$ be a polygonal open bounded subset of $\R^2$, and
let $\partial_D \Om \subseteq \partial \Om$ be open in the relative topology.
For all $\varepsilon>0$, let $\rb_\varepsilon \in \rs_\varepsilon(\Om)$ be
a regular triangulation of $\Om$ such that $\partial_D \Om$ is composed of edges of $\rb_\eps$.
As in the previous section, let $\afeom$ be the family of continuous piecewise affine
functions with respect to $\rb_\eps$, and let $\afeaom$ be the family of functions which are affine
on the triangles of some triangulation $\tb \in \treaom$ nested in $\rb_\eps$ and can jump across
the edges of $\tb$.
\par
In the following, it will be useful to treat points at which the boundary condition
is violated (see \eqref{jump*}) as internal jumps. Thus we
consider $\Om_D$ polygonal open bounded subset of $\R^2$ such that
$\Om_D \cap \Om=\emptyset$ and $\partial \Om \cap \partial \Om_D= \partial_D \Om$ up to a finite number of
points;
we set $\Om':= \Om \cup \Om_D \cup \partial_D \Om$.
Given $u \in \afeaom$ and $g \in \afeom$, we may extend
$g$ to a function of $H^1(\Om')$ and $u$ to a function
$\tilde{u} \in SBV(\Om')$ setting $\tilde{u}=g$ on $\Om_D$.
In this way, recalling \eqref{jumpbis}, we have
\begin{equation*}
\Sg{g}{u}=S_{\tilde{u}},
\end{equation*}
so that the violation of the boundary condition of $u$ can be read in the
set of jumps of $\tilde{u}$.
Analogously, given $u \in SBV(\Om)$ and $g \in H^1(\Om)$, we set
\begin{equation}
\label{jumptris}
\Sg{g}{u}:=S_u \cup \{x \in \partial_D \Om\,:\, \gamma(u)(x) \not= \gamma(g)(x)\}
\end{equation}
where $\gamma$ denotes the trace operator on $\partial \Om$.
We may assume $g \in H^1(\Om')$ using an extension operator.
We can then consider $\tilde{u} \in SBV(\Om')$
such that $\tilde{u}=u$ on $\Om$, and $\tilde{u}=g$ on $\Om_D$.
In this way we have
\begin{equation*}
\Sg{g}{u}=S_{\tilde{u}} \quad \mbox{ up to a set of $\hs^1${-}measure $0$}.
\end{equation*}
Let us consider $g \in W^{1,1}([0,1],H^1(\Om))$ such that $\|g(t)\|_\infty \le C$
for all $t \in [0,1]$ and let $g_\eps \in W^{1,1}([0,1],H^1(\Om))$ be such that
$g_\eps(t) \in \afeom$ for all $t \in [0,1]$,
\begin{equation}
\label{lib}
\|g_\eps(t)\|_\infty \le C
\end{equation}
for all $t \in [0,1]$, and for $\eps \to 0$
\begin{equation}
\label{strconv}
g_\eps \to g
\quad \mbox{ strongly in }W^{1,1}([0,1],H^1(\Om)).
\end{equation}
We indicate by $\{u^{\delta,i}_{\eps,a},\,i=0, \ldots,N_\delta\}$ the
discrete evolution relative to the boundary data $g_\eps$ given by Proposition
\ref{discrevol}, and we denote by $\Es^{\delta,i}_{\eps,a}$ its total energy as in
Proposition \ref{dener}.
\par
We assume that $g(\cdot)$ and $g_\eps(\cdot)$ are defined in
$H^1(\Om')$ (we still denote these extensions by $g(\cdot)$ and $g_h(\cdot)$), in
such a way that \eqref{lib} and \eqref{strconv} hold in $\Om'$.
Let us moreover set $g^\delta_\eps(t):=g_\eps(t^\delta_i)$
for all $t^\delta_i \le t <t^\delta_{i+1}$ with $i=0,\dots,N_\delta-1$ and
$g^\delta_\eps(1):=g_\eps(1)$.
\par
Let us make the following piecewise constant interpolation in time:
\begin{equation*}
u^\delta_{\eps,a}(t):=u^{\delta,i}_{\eps,a}
\quad \mbox{ for } t^\delta_i \le t <t^\delta_{i+1}
\quad i=0,\dots,N_\delta-1,
\end{equation*}
and $u^\delta_{\eps,a}(1):=u^{\delta,N_\delta}_{\eps,a}$.
For all $t \in [0,1]$ we define the {\it discrete fracture} at time $t$ as
\begin{equation*}
\Gamma^\delta_{\eps,a}(t):= \bigcup_{s \le t}
\Sg{g^\delta_\eps(s)}{u^\delta_{\eps,a}(s)},
\end{equation*}
and the {\it discrete total energy} at time $t$ as
\begin{equation*}
\Es^\delta_{\eps,a}(t):=
\|\nabla u^\delta_{\eps,a}(t)\|^2+
\hs^1 \left( \Gamma^\delta_{\eps,a}(t) \right).
\end{equation*}
We have for all $t \in [0,1]$
\begin{equation}
\label{boundonu}
\|u^{\delta}_{\eps,a}(t)\|_\infty \le \|g_\eps^\delta(t)\|_\infty.
\end{equation}
Moreover for all $v \in \afeaom$ we have
\begin{equation}
\label{MSdiscr2}
\|\nabla u^{\delta}_{\eps,a}(0)\|^2
+\hs^1\left( \Sg{g^\delta_\eps(0)}{u^{\delta}_{\eps,a}(0)} \right)
\le \|\nabla v\|^2 +
\hs^1\left( \Sg{g^\delta_\eps(0)}{v} \right),
\end{equation}
and for all $t \in ]0,1]$ and for all $v \in \afeaom$
\begin{equation}
\label{piecemin2}
\|\nabla u^{\delta}_{\eps,a}(t)\|^2 \le \|\nabla v\|^2 +
\hs^1 \left( \Sg{g^\delta_\eps(t)}{v}  \setminus
\Gamma^{\delta}_{\eps,a}(t) \right).
\end{equation}
Finally for all $0 \le s \le t \le 1$ we have
\begin{equation}
\label{discrenergybists}
\Es^{\delta}_{\eps,a}(t) \le \Es^{\delta}_{\eps,a}(s)+
2\int_{s^\delta_{i}}^{t^\delta_{i}}
\int_{\Om} \nabla u^{\delta}_{\eps,a}(\tau) \nabla
\dot{g}_\eps(\tau) \,dx\,d\tau +o^\delta_{\eps},
\end{equation}
where $t^\delta_i \le t <t^\delta_{i+1}$,
$s^\delta_i \le s <s^\delta_{i+1}$ and
\begin{equation}
\label{opiccolo}
o^\delta_{\eps}:=
\left[
\max_{r=0,\ldots,N_\delta-1}
\int_{t^\delta_{r}}^{t^\delta_{r+1}} \|\dot{g}_\eps(\tau)\|_{H^1(\Om)} \,d\tau
\right] \int_0^1 \|\dot{g}_\eps(\tau)\|_{H^1(\Om)}.
\end{equation}
For $s=0$ we obtain the following estimate from above for the discrete total energy
\begin{equation}
\label{discrenergybis}
\Es^{\delta}_{\eps,a}(t) \le \Es^{\delta}_{\eps,a}(0)+
2\int_{0}^{t^\delta_{i}}
\int_{\Om} \nabla u^{\delta}_{\eps,a}(\tau) \nabla
\dot{g}_\eps(\tau) \,dx\,d\tau +o^\delta_{\eps},
\end{equation}
where $t^\delta_i \le t <t^\delta_{i+1}$.
\par
We study the behavior of the evolution
$\{t \to u^\delta_{\varepsilon,a}(t),\,t \in [0,1]\}$ varying the parameters
in the following way. We let firstly $\eps \to 0$ and $\delta \to 0$
obtaining an evolution $\{t \to u_a(t),\,t \in [0,1]\}$
relative to the boundary data $g$ with the minimality property \eqref{amint}; then we
let $a \to 0$ obtaining a quasi-static evolution of brittle fractures
$\{t \to u(t)\,,t \in [0,1]\}$ relative to the boundary data $g$.
Finally, by a diagonal argument we deal with $(\delta,\eps,a)$ at the
same time.
\par
In order to develop this program, we need some compactness, and so
we derive a bound for the total energy $\Es^\delta_{\eps,a}$.
By \eqref{piecemin}, we have that for all $t \in [0,1]$
\begin{equation*}
\|\nabla u^{\delta}_{\varepsilon,a}(t)\| \le \|\nabla g^\delta_\eps(t)\| \le \tilde{C}
\end{equation*}
with $\tilde C$ independent of $\delta$, $\eps$ and $t$.
We deduce for all $t \in [0,1]$
\begin{equation*}
\Es^{\delta}_{\varepsilon,a}(t) \le
\Es^{\delta}_{\varepsilon,a}(0)+
2 \tilde{C}^2 +o^\delta_\eps
\end{equation*}
Notice that $\Es^\delta_{\eps,a}(0)$ is uniformly bounded as
$\delta,\eps$ vary. Moreover, by \eqref{boundonu} and since
$\|g_\eps(t)\|_\infty \le C$ for all $t \in [0,1]$, we have that
$u^{\delta}_{\varepsilon,a}(t)$
is uniformly bounded in $L^\infty(\Om)$ independently of $\delta,\eps$ and $a$.
Taking into account \eqref{strconv}, we conclude
that there exists $C'$ independent of $\delta,\varepsilon,a$  such that for all
$t\in [0,1]$
\begin{equation}
\label{unifenergybound}
\Es^{\delta}_{\varepsilon,a}(t) + \|u^{\delta}_{\varepsilon,a}(t)
\|_\infty \le C'.
\end{equation}
Formula \eqref{unifenergybound} gives the desired compactness in order to perform
the asymptotic analysis of the discrete evolution.
\par
Let now consider
$\delta_n \to 0$ and $\eps_n \to 0$: by \eqref{strconv} we have
\begin{equation}
\label{Cto0}
o^{\delta_n}_{\eps_n} \to 0,
\end{equation}
where $o^{\delta_n}_{\eps_n}$ is defined in \eqref{opiccolo}.
By Helly's theorem on monotone functions, we may suppose that there exists an
increasing function $\lambda_a$ such that (up
to a subsequence) for all $t \in [0,1]$
\begin{equation}
\label{lambn}
\lambda_{n,a}(t):=
\hs^1 \left( \bigcup_{s \le t}\Sg{g^{\delta_n}_{\eps_n}(s)}
{u^{\delta_n}_{\eps_n,a}(s)} \right)
\to \lambda_a(t).
\end{equation}
Let us fix $D \subseteq [0,1]$ countable and dense with $0 \in D$.

\begin{lemma}
\label{eps-to-0}
For all $t \in D$ there exists $\vuat \in SBV(\Om)$
such that up to a subsequence independent of $t$
\begin{equation*}
u^{\delta_n}_{\eps_n,a}(t) \to \vuat \quad \mbox{ in } SBV(\Om).
\end{equation*}
Moreover for all $t \in D$ we have
\begin{equation}
\label{ueb2}
\|\nabla \vuat \|^2+ \hs^1 \big( \Sg{g(t)}{\vuat} \big) +
\|\vuat\|_\infty \le C'.
\end{equation}
\end{lemma}

\begin{proof}
Let us consider $t \in D$. By \eqref{unifenergybound}, we can apply
Ambrosio's Compactness Theorem \ref{SBVcompact} obtaining $u \in SBV(\Om)$ such
that, up to a subsequence, $u^{\delta_n}_{\eps_n,a}(t) \to u$ in $SBV(\Om)$.
Let us set $\vuat:=u$. Using a diagonal argument, we deduce that there exists a
subsequence of $(\delta_n, \eps_n)$ (which we still denote by
$(\delta_n,\eps_n)$) such that $u^{\delta_n}_{\eps_n,a}(t) \to
\vuat$ in $SBV(\Om)$ for all $t \in D$.
In order to obtain inequality \eqref{ueb2}, we extend
$u^{\delta_n}_{\eps_n,a}(t)$ and $\vuat$ to $\Om'$ setting
$u^{\delta_n}_{\eps_n,a}(t):= g^{\delta_n}_{\eps_n}(t)$ and
$u_a(t):=g(t)$ on $\Om_D$; since $g^{\delta_n}_{\eps_n}(t) \to g(t)$ on
$\Om_D$ strongly in $H^1(\Om_D)$, we have that
$u^{\delta_n}_{\eps_n,a}(t) \to u_a(t)$ in $SBV(\Om')$, so that
we can apply Ambrosio's Theorem, and derive \eqref{ueb2} from
\eqref{unifenergybound}.
\end{proof}

The following result is essential for the sequel: its proof is postponed to
Section \ref{secmin}.

\begin{proposition}
\label{pminpropD}
Let $t \in D$. For all $v \in SBV(\Om)$ we have
\begin{equation}
\label{minpropD}
\|\nabla \vuat\|^2 \le \|\nabla v\|^2+ \mu(a) \hs^1(\Sg{g(t)}{v} \setminus
\bigcup_{s \le t, s \in D} \Sg{g(s)}{\vuas}),
\end{equation}
where $\mu:]0,\frac{1}{2}[ \to ]0,+\infty[$ is such that
$\lim_{a \to 0} \mu(a)=1$.
Moreover, $\nabla u^{\delta_n}_{\eps_n,a}(t) \to \nabla \vuat$ strongly in
$L^2(\Om;\R^2)$.
\end{proposition}

We now extend the evolution $\{t \to u_a(t)\,:\,t \in D\}$ to the entire interval
$[0,1]$. Let us set for all $t \in [0,1]$
\begin{equation*}
\Gamma_a(t):= \bigcup_{s \le t, s \in D} \Sg{g(s)}{u_a(s)}.
\end{equation*}

\begin{lemma}
\label{aextension}
For every $t \in [0,1]$ there exists $u_a(t) \in SBV(\Om)$
such that the following hold:
\begin{itemize}
\item[(a)] for all $t \in [0,1]$
\begin{equation}
\label{ajumpt}
\Sg{g(t)}{u_a(t)} \subseteq
\Gamma_a(t) \;\mbox{ up to a set of }\hs^1\mbox{-measure }0,
\end{equation}
and
\begin{equation}
\label{aueb}
\|\nabla \vuat \|^2+ \hs^1 \big( \Sg{g(t)}{\vuat} \big) +
\|\vuat\|_\infty \le C';
\end{equation}
\item[{}]
\item[(b)] for all $v \in SBV(\Om)$
\begin{equation}
\label{amint}
\|\nabla u_a(t)\|^2 \le \|\nabla v\|^2 +
\mu(a)\hs^1 \left( \Sg{g(t)}{v} \setminus \Gamma_a(t) \right);
\end{equation}
\item[{}]
\item[(c)]
$\nabla u_a$ is left continuous in $[0,1] \setminus D$ with respect to the
strong topology of $L^2(\Om;\R^2)$;
\item[{}]
\item[(d)]
for all $t \in [0,1] \setminus \ns_a$ we have that
\begin{equation*}
\nabla u^{\delta_n}_{\eps_n,a}(t) \to \nabla u_a(t)
\quad \mbox{ strongly in }L^2(\Om, \R^2),
\end{equation*}
where $\ns_a$ is the set of discontinuities of the function $\lambda_a$
defined in \eqref{lambn}.
\end{itemize}
\end{lemma}

\begin{proof}
Let $t \in [0,1] \setminus D$ and let $t_n \in D $ with $t_n \nearrow t$.
By \eqref{ueb2}, we can apply Ambrosio's Theorem to the sequence $(u_a(t_n))$
obtaining $u \in SBV(\Om)$ such that, up to a subsequence, $u_a(t_n) \to u$ in
$SBV(\Om)$. Let us set $\vuat:=u$.
Let us extend $u_a(t_n)$ and $\vuat$ to $\Om'$ setting $u_a(t_n):=g(t_n)$ and
$\vuat:=g(t)$ on $\Om_D$: we have $u_a(t_n) \to \vuat$ in $SBV(\Om')$.
Since $\hs^1 \res S_{u_a(t_n)} \le \hs^1 \res \Gamma_a(t)$ for all $n$,
as a consequence of Ambrosio's Theorem, we deduce that
$\hs^1 \res S_{\vuat} \le \hs^1 \res \Gamma_a(t)$. This means
$\hs^1 \res \Sg{g(t)}{\vuat} \le \hs^1 \res \Gamma_a(t)$, so that
\eqref{ajumpt} holds.
Moreover, for all $v \in SBV(\Om)$, by \eqref{minpropD} we may write
\begin{multline}
\label{minbis}
\|\nabla u_a(t_n)\|^2 \le
\|\nabla v-\nabla g(t)+\nabla g(t_n)\|^2+
\mu(a) \hs^1 \left( \Sg{g(t)}{v} \setminus \Gamma_a(t_n) \right) \le \\
\le \|\nabla v-\nabla g(t)+\nabla g(t_n)\|^2+
\mu(a) \hs^1 \left( \Sg{g(t)}{v} \setminus \Gamma_a(t) \right) +
\mu(a) \hs^1 \left( \Gamma_a(t) \setminus \Gamma_a(t_n) \right),
\end{multline}
so that, since by definition of $\Gamma_a(t)$ we have
$\hs^1(\Gamma_a(t) \setminus \Gamma_a(t_n)) \to 0$, we obtain that
\eqref{amint} holds; choosing $v=u_a(t)$ and taking the limsup in \eqref{minbis},
we obtain that
$$
\limsup_n \|\nabla u_{a}(t_n)\|^2 \le \|\nabla u_a(t)\|^2,
$$
and so the convergence $\nabla u_a(t_{n}) \to \nabla u_a(t)$ is strong in
$L^2(\Om, \R^2)$.
Notice that $\nabla u_a(t)$ is uniquely determined by (\ref{ajumpt}) and (\ref{amint})
since the gradient of the solutions of the minimum problem
\begin{equation*}
\min \left\{ \|\nabla u\|^2 \,:\, \Sg{g(t)}{u} \subseteq \Gamma_a(t)
\mbox{ up to a set of }\hs^1\mbox{-measure }0 \right\}
\end{equation*}
is unique by the strict convexity of the functional: we conclude that $\nabla u_a(t)$
is well defined.
The same arguments prove that $\nabla u_a$ is left continuous
at all $t \in [0,1] \setminus D$.
Finally \eqref{aueb} is a direct
consequence of \eqref{ueb2} and of Ambrosio's Theorem, and so points $(a)$, $(b)$, $(c)$ are
proved.
\par
Let us come to point $(d)$.
Let us consider $u^{\delta_n}_{\eps_n,a}(t)$ with $t \not\in \ns_a$; we
may suppose that $t \not \in D$, since otherwise the result has already been
established. By Proposition \ref{pminpropD} with $D':=D \cup\{t\}$ in place of $D$,
we have that, up to
a subsequence, $u^{\delta_n}_{\eps_n,a}(t) \to u$ in $SBV(\Om)$
such that
\begin{equation*}
\|\nabla u\|^2 \le \|\nabla v\|^2 +\mu(a) \hs^1 \left( \Sg{g(t)}{v}
\setminus (\Gamma_a(t) \cup \Sg{g(t)}{u}) \right)
\end{equation*}
for all $v \in SBV(\Om)$ and
$\nabla u^{\delta_n}_{\eps_n,a}(t) \to \nabla u$ strongly in  $L^2(\Om;\R^2)$.
Let $s<t$ with $ s\in D$; by the minimality of $u^{\delta_n}_{\eps_n,a}(s)$
and by \eqref{unifenergybound} we have
\begin{multline*}
\|\nabla u^{\delta_n}_{\eps_n,a}(s)\|^2 \le
\|\nabla u^{\delta_n}_{\eps_n,a}(t)-\nabla g^{\delta_n}_{\eps_n}(t)+
\nabla g^{\delta_n}_{\eps_n}(s)\|^2
+\lambda_{n,a}(t)-\lambda_{n,a}(s) \le \\
\le \|\nabla u^{\delta_n}_{\eps_n,a}(t)\|^2
+2 \sqrt{C'}\|\nabla g^{\delta_n}_{\eps_n}(t)-\nabla g^{\delta_n}_{\eps_n}(s)\|+ \\
+\|\nabla g^{\delta_n}_{\eps_n}(t)-\nabla g^{\delta_n}_{\eps_n}(s)\|^2+
\lambda_{n,a}(t)-\lambda_{n,a}(s).
\end{multline*}
Passing to the limit for $n \to +\infty$, recalling that
$g^{\delta_n}_{\eps_n}(\tau) \to g(\tau)$
strongly in $H^1(\Om)$ for all $\tau \in [0,1]$, we deduce
\begin{multline*}
\|\nabla u_a(s)\|^2 \le \|\nabla u\|^2+
2\sqrt{C'}\|\nabla g(t)-\nabla g(s)\|+
\|\nabla g(t)-\nabla g(s)\|^2+
\lambda_a(t)-\lambda_a(s),
\end{multline*}
so that, since $t$ is a point of continuity for $\lambda_a$, $\nabla u_a$
is left continuous at $t$, and $g$ is absolutely continuous, we get
for $s \to t$
\begin{equation*}
\|\nabla u_a(t)\|^2 \le \|\nabla u\|^2.
\end{equation*}
We conclude that $u_a(t)$ is a solution of
\begin{equation*}
\min \{\|\nabla v\|^2\,:\,
\Sg{g(t)}{v} \subseteq \Gamma_a(t) \cup \Sg{g(t)}{u}
\mbox{ up to a set of }\hs^1\mbox{-measure }0\},
\end{equation*}
so that $\nabla u=\nabla u_a(t)$ by uniqueness of the gradient of the solution.
We deduce that $\nabla u^{\delta_n}_{\eps_n,a}(t) \to \nabla u_a(t)$
strongly in $L^2(\Om;\R^2)$, and so the proof is complete.
\end{proof}

We can now let $a \to 0$.

\begin{lemma}
\label{azero}
There exists $a_n \to 0$ such that, for all $t \in D$,
$u_{a_n}(t) \to u(t)$ in $SBV(\Om)$ for some
$u(t) \in SBV(\Om)$ such that for all $v \in SBV(\Om)$
we have
\begin{equation}
\label{alessminpropD}
\|\nabla u(t)\|^2 \le \|\nabla v\|^2+ \hs^1(\Sg{g(t)}{v} \setminus
\bigcup_{s \le t, s \in D} \Sg{g(s)}{u(s)}).
\end{equation}
Moreover, $\nabla u_{a_n}(t) \to \nabla u(t)$ strongly in $L^2(\Om;\R^2)$ and
\begin{equation}
\label{ueb3}
\|\nabla u(t) \|^2+ \hs^1(\Sg{g(t)}{u(t)}) +
\|u(t)\|_\infty \le C'.
\end{equation}
\end{lemma}

\begin{proof}
By \eqref{aueb}, applying Ambrosio's Theorem to the extensions of
$\vuat$ to $\Om'$ by setting $\vuat:=g(t)$ on $\Om_D$, and using a diagonal argument,
we find a sequence $a_n \to 0$ such that, for all $t \in D$,
$u_{a_n}(t) \to u(t)$ in $SBV(\Om)$ for some
$u(t) \in SBV(\Om)$ such that \eqref{ueb3} holds.
\par
We now prove that $u(t)$ satisfies property (\ref{alessminpropD}).
Let $v \in SBV(\Om)$.
Let us fix $t_1 \le t_2 \le \ldots \le t_k=t$ with $t_i \in D$.
We extend $v$ and $u_{a_n}(t_i)$ to $\Om'$ setting $v:=g(t)$ and
$u_{a_n}(t_i):=g(t_i)$ on $\Om_D$ respectively.
Since $u_{a_n}(t_i) \to u(t_i)$ in $SBV(\Om')$ for all $i=1, \ldots, k$, by
Theorem \ref{jumptransfer} there exists $v_n \in SBV(\Om')$ with $v_n=g(t)$ on
$\Om_D$ such that $\nabla v_n \to \nabla v$ strongly in $L^2(\Om';\R^2)$ and
\begin{equation}
\label{convlength}
\limsup_n \hs^1 \left( S_{v_n} \setminus \bigcup_{i=1}^k S_{u_{a_n}(t_i)} \right)
\le \hs^1 \left( S_v \setminus \bigcup_{i=1}^k S_{u(t_i)} \right).
\end{equation}
By \eqref{minpropD} we obtain
\begin{equation}
\label{minan}
\|\nabla u_{a_n}(t)\|^2 \le \|\nabla v_n\|^2+ \mu(a_n)
\hs^1 \left( S_{v_n} \setminus \bigcup_{i=1}^k S_{u_{a_n}(t_i)} \right),
\end{equation}
so that passing to the limit for $n \to +\infty$ and recalling that
$\mu(a) \to 1$ as $a \to 0$, we obtain
$$
\|\nabla u(t)\|^2 \le \|\nabla v\|^2+
\hs^1 \left( S_v \setminus \bigcup_{i=1}^k S_{u(t_i)} \right).
$$
Thus we get
$$
\|\nabla u(t)\|^2 \le \|\nabla v\|^2+
\hs^1 \left( \Sg{g(t)}{v} \setminus \bigcup_{i=1}^k \Sg{g(t_i)}{u(t_i)} \right).
$$
Since $t_1, \ldots, t_k$ are arbitrary, we obtain \eqref{alessminpropD}.
Choosing $v=u(t)$, taking the limsup in \eqref{minan} and using
\eqref{convlength}, we obtain $\nabla u_{a_n}(t) \to \nabla u(t)$ strongly in
$L^2(\Om;\R^2)$.
\end{proof}

In order to deal with $\delta,\eps$ and $a$ at the same time,
we need the following lemma.

\begin{lemma}
\label{convD}
Let $\{u(t)\,:\,t \in D\}$ be as in Lemma \ref{azero}.
There exist $\delta_n \to 0$, $\varepsilon_n \to 0$, and $a_n \to 0$
such for all $t \in D$ we have
\begin{equation*}
u^{\delta_n}_{\varepsilon_n,a_n}(t) \to u(t)
\quad \mbox{ in $SBV(\Om)$ }.
\end{equation*}
Moreover, for all $n$ there exists $\bs_n \subseteq [0,1]$
with $|\bs_n|<2^{-n}$ such that for all $t \in [0,1] \setminus \bs_n$
\begin{equation}
\label{graduanconv}
\|\nabla u^{\delta_n}_{\varepsilon_n,a_n}(t)-\nabla u_{a_n}(t)\|
\le \frac{1}{n}.
\end{equation}
Finally, we have that for all $v \in SBV(\Om)$
\begin{equation}
\label{min0}
\| \nabla u(0)\|^2 +\hs^1 \left( \Sg{g(0)}{u(0)} \right)
\le \|\nabla v\|^2+\hs^1 \left( \Sg{g(0)}{v} \right)
\end{equation}
and
\begin{equation}
\label{conv0}
\Es^{\delta_n}_{\eps_n,a_n}(0) \to \|\nabla u(0)\|^2+
\hs^1 \left( \Sg{g(0)}{u(0)} \right).
\end{equation}
\end{lemma}

\begin{proof}
Let $(a_n)$ be the sequence determined by Lemma \ref{azero}.
By Lemma \ref{eps-to-0}, for all $n$ there exists
$(\delta^n_m,\varepsilon^n_{m})$
such that for all $t \in D$ and $m \to +\infty$ we have
\begin{equation*}
u^{\delta^n_m}_{\varepsilon^n_m,a_n}(t) \to u_{a_n}(t)
\quad \mbox{ in }SBV(\Om),
\end{equation*}
and
\begin{equation*}
\nabla u^{\delta^n_m}_{\varepsilon^n_m,a_n}(t)
\to \nabla u_{a_n}(t)
\quad \mbox{ strongly in }L^2(\Om;\R^2).
\end{equation*}
Moreover by Lemma \ref{aextension} we have that
$\nabla u^{\delta^n_m}_{\varepsilon^n_m,a_n} \to \nabla u_{a_n}$
quasi-uniformly on $[0,1]$ as $m \to +\infty$.
Let $\bs_n \subseteq [0,1]$ with $|\bs_n| <2^{-n}$ such that
$\nabla u^{\delta^n_m}_{\varepsilon^n_m,a_n} \to \nabla u_{a_n}$
uniformly on $[0,1] \setminus \bs_n$ as $m \to +\infty$.
We now perform the following diagonal argument. Let $D=\{t_n,\,n \ge 1\}$.
Choose $m_1$ such that
\begin{equation*}
\|\nabla u^{\delta^1_{m_1}}_{\varepsilon^1_{m_1},a_1}(t_1)-
\nabla u_{a_1}(t_1)\|+
\|u^{\delta^1_{m_1}}_{\varepsilon^1_{m_1},a_1}(t_1)-u_{a_1}(t_1)\|
\le 1,
\end{equation*}
and
\begin{equation*}
\|\nabla u^{\delta^1_{m_1}}_{\varepsilon^1_{m_1},a_1}(t)-
\nabla u_{a_1}(t)\| \le 1
\quad \mbox{ for all $t \in [0,1] \setminus \bs_1$ }.
\end{equation*}
Let $m_n$ be such that
\begin{equation*}
\|\nabla u^{\delta^n_{m_n}}_{\varepsilon^n_{m_n},a_n}(t_j)
-\nabla u_{a_n}(t_j)\|+
\|u^{\delta^n_{m_n}}_{\varepsilon^n_{m_n},a_n}(t_j)-u_{a_n}(t_j)\|
\le \frac{1}{n}
\;\;\mbox{ for all $j=1, \ldots, n$}
\end{equation*}
and
\begin{equation*}
\|\nabla u^{\delta^n_{m_n}}_{\varepsilon^n_{m_n},a_n}(t)-
\nabla u_{a_n}(t)\| \le \frac{1}{n}
\quad \mbox{ for all $t \in [0,1] \setminus \bs_n$ }.
\end{equation*}
We may suppose that $\delta^n_{m_n} \to 0$, $\eps^n_{m_n} \to 0$.
Then $(\delta^n_{m_n},\eps^n_{m_n},a_n)$ is the sequence which satisfies
the thesis. In fact by construction and taking into account \eqref{unifenergybound},
for all $t \in D$ we have
$u^{\delta^n_{m_n}}_{\varepsilon^n_{m_n},a_n}(t) \to u(t)$ in $SBV(\Om)$;
moreover the set $\bs_n$ satisfies \eqref{graduanconv}.
Notice that $u^{\delta^n_{m_n}}_{\eps^n_{m_n},a^n_{m_n}}(0)$
satisfies \eqref{MSdiscr2} and so \eqref{min0} and \eqref{conv0} follow by
the $\Gamma${-}convergence result of \cite{N}.
\end{proof}

Let $(\delta_n,\eps_n,a_n)$ be the sequence determined by Lemma
\ref{convD}. For all $t \in [0,1]$ let us set
\begin{equation*}
\lambda_n(t):=
\hs^1 \left( \Gamma^{\delta_n}_{\eps_n,a_n}(t) \right).
\end{equation*}
By Helly's theorem, we may suppose that there exist two increasing functions
$\lambda$ and $\eta$ such that up to a subsequence
\begin{equation*}
\lambda_n \to \lambda \quad \mbox{ pointwise in }[0,1],
\end{equation*}
and
\begin{equation}
\label{eta}
\lambda_{a_n} \to \eta \quad \mbox{ pointwise in }[0,1],
\end{equation}
where $\lambda_{a_n}$ is defined as in \eqref{lambn}.
We now extend the evolution $\{t \to u(t)\,:\,t \in D\}$ to the entire interval
$[0,1]$. Let us set for all $t \in [0,1]$
\begin{equation*}
\Gamma(t):= \bigcup_{s \le t, s \in D} \Sg{g(s)}{u(s)},
\end{equation*}
and let $\ns$ be the set of discontinuities of $\hs^1(\Gamma(\cdot))$.
Notice that for all $t \in [0,1]$
\begin{equation}
\label{size2}
\hs^1(\Gamma(t)) \le \lambda(t).
\end{equation}
In fact if $t \in D$, let $t_1 \le t_2 \le \ldots \le t_k=t$ with $t_i \in D$,
consider $w_n \in SBV(\Om';\R^k)$ defined as
\begin{equation*}
w_n(x):=(u^{\delta_n}_{\varepsilon_n,a_n}(t_1)(x), \ldots,
u^{\delta_n}_{\varepsilon_n,a_n}(t_k)(x)),
\end{equation*}
where we assume that $u^{\delta_n}_{\varepsilon_n,a_n}(t_i)=g^{\delta_n}_{\eps_n}(t_i)$
on $\Om_D$. We have $w_n \to w:=(u(t_1), \ldots,  u(t_k))$ in $SBV(\Om';\R^k)$,
where $u(t_i)=g(t_i)$ on $\Om_D$.
Note that for all $n$ we have $S_{w_n}=
\bigcup_{i=1}^k S_{u^{\delta_n}_{\varepsilon_n,a_n}(t_i)}$ so that
\begin{equation*}
\hs^1(S_{w_n}) \le \lambda_n(t).
\end{equation*}
Passing to the limit for $n \to +\infty$ and applying Ambrosio's Theorem we get
\begin{equation*}
\hs^1 \left( \bigcup_{i=1}^k S_{u(t_i)} \right)
=\hs^1(S_w) \le \liminf_n \hs^1(S_{w_n}) \le
\lambda(t);
\end{equation*}
we thus have
\begin{equation*}
\hs^1 \left( \bigcup_{i=1}^k \Sg{g(t_i)}{u(t_i)} \right)
=\hs^1(S_w) \le \lambda(t)
\end{equation*}
and taking the sup over all $t_1, \ldots, t_k$, we obtain \eqref{size2} in $D$.
The case $t \not \in D$ follows since $\hs^1(\Gamma(\cdot))$ is left continuous
by definition.

\begin{lemma}
\label{extension}
For every $t \in [0,1]$ there exists $u(t) \in SBV(\Om)$
such that the following hold:
\begin{itemize}
\item[(a)] for all $t \in [0,1]$
\begin{equation}
\label{jumpt}
\Sg{g(t)}{u(t)} \subseteq
\Gamma(t) \; \mbox{ up to a set of }\hs^1\mbox{-measure }0,
\end{equation}
and for all $t \in [0,1]$ and for all $v \in SBV(\Om)$
\begin{equation}
\label{mint}
\|\nabla u(t)\|^2 \le
\|\nabla v\|^2+ \hs^1 \left( \Sg{g(t)}{v} \setminus \Gamma(t) \right);
\end{equation}
\item[{}]
\item[(b)]  $\nabla u$ is continuous in $[0,1] \setminus (D \cup \ns)$ with respect to
the strong topology of $L^2(\Om;\R^2)$;
\item[{}]
\item[(c)] if $\tilde{\ns}$ is the set of discontinuities of the function $\eta$ defined in
\eqref{eta}, for all $t \in [0,1] \setminus \tilde{\ns}$
we have that
\begin{equation*}
\nabla u_{a_n}(t) \to \nabla u(t)
\quad \mbox{ strongly in }L^2(\Om, \R^2).
\end{equation*}
\end{itemize}
Finally
\begin{equation}
\label{energybelow}
\Es(t) \ge \Es(0)+
2\int_0^t \int_{\Om} \nabla u(\tau) \nabla \dot{g}(\tau) \,dx\,d\tau,
\end{equation}
where
\begin{equation*}
\label{defenergy}
\Es(t):= \|\nabla u(t)\|^2+\hs^1(\Gamma(t)).
\end{equation*}
\end{lemma}

\begin{proof}
The definition of $u(t)$ is carried out as in Lemma \ref{aextension}
considering $t \in [0,1] \setminus D$, $t_n \in D$ with $t_n \nearrow t$,
and the limit (up to a subsequence) of $u(t_n)$ in $SBV(\Om)$:
\eqref{jumpt} and \eqref{mint} hold, so that point $(a)$ is proved.
It turns out that $\nabla u(t)$ is uniquely determined and that it is left continuous
in $[0,1] \setminus D$.
Let us consider $t \in [0,1] \setminus (D \cup \ns)$, and let $t_n \searrow t$.
By Ambrosio's Theorem, we have that there exists $u \in SBV(\Om)$ with
such that, up to a subsequence, $u(t_n) \to u$ in $SBV(\Om)$.
Since $t$ is a continuity point of $\hs^1(\Gamma(\cdot))$, we deduce that
$\Sg{g(t)}{u} \subseteq \Gamma(t)$ up to a set of $\hs^1$-measure $0$.
Moreover by the minimality property for $u(t_n)$ and the fact
$\Gamma(t) \subseteq \Gamma(t_n)$, we have that for all $v \in SBV(\Om)$ with
\begin{multline*}
\|\nabla u(t_n)\|^2 \le
\|\nabla v-\nabla g(t)+\nabla g(t_n)\|^2+
\hs^1 \left( \Sg{g(t)}{v} \setminus \Gamma(t_n) \right) \le \\
\le \|\nabla v-\nabla g(t)+\nabla g(t_n)\|^2+
\hs^1 \left( \Sg{g(t)}{v} \setminus \Gamma(t) \right),
\end{multline*}
and so we deduce that \eqref{mint} holds with $u$ in place of $u(t)$, and that
$\nabla u(t_n) \to \nabla u$ strongly in $L^2(\Om;\R^2)$.
We obtain by uniqueness that $\nabla u=\nabla u(t)$, and so $\nabla u(\cdot)$ is continuous in
$[0,1] \setminus (D \cup \ns)$ and this proves point $(b)$.
Point $(c)$ follows in the same way of point $(d)$ of Lemma \ref{aextension}.
\par
Let us come to the proof of \eqref{energybelow}. Given $t \in [0,1]$ and $k>0$,
let $s_i^k:= \frac{i}{k}t$ for all $i=0, \ldots, k$. Let us set
$u^k(s):=u(s_{i+1}^k)$ for $s_i^k< s \le s_{i+1}^k$. By (\ref{mint}), comparing
$u(s^k_i)$ with $u(s^k_{i+1})-g(s^k_{i+1})+g(s^k_{i})$, it is easy to see that
\begin{equation*}
\Es(t) \ge \Es(0)+
2\int_0^t \int_{\Om} \nabla u^k(\tau) \nabla \dot{g}(\tau) \,d\tau \,dx+o_k,
\end{equation*}
where $o_k \to 0$ as $k \to +\infty$.
Since $\nabla u$ is continuous with respect to the strong topology of $L^2(\Om;\R^2)$
in $[0,1]$ up to a countable set, passing to the limit for $k \to +\infty$ we deduce
\eqref{energybelow}.
\end{proof}

We are now ready to prove the main result of the paper.
\begin{proof}[\rm \underline{PROOF OF THEOREM \ref{mainthm}}]
Let $D$ be a countable and dense set in $[0,1]$ such that
$0 \in D$, and let $(\delta_n,\eps_n,a_n)$ and
$\{t \to u(t) \in SBV(\Om)\,:\, t \in [0,1]\}$
be the sequence and the evolution determined in Lemma \ref{convD} and
Lemma \ref{extension}.
Let us set
$$
u_n:=u^{\delta_n}_{\varepsilon_n, a_n}, \quad \quad
\Gamma_n:=\Gamma^{\delta_n}_{\varepsilon_n, a_n}, \quad \quad
\Es_n:=\Es^{\delta_n}_{\varepsilon_n, a_n}.
$$
Let $\overline{\ns}$ be the union of the sets of discontinuities of
$\eta$ and $\hs^1(\Gamma(\cdot))$,
where $\eta$ is defined in \eqref{eta}.
Let $\bs:= \bigcap_{k=1}^{+\infty} \bigcup_{h=k}^\infty \bs_h$, where
$\bs_h$ are as in Lemma \ref{convD}; since
$|\bigcup_{h=k}^\infty \bs_h| <2^{-k+1}$, we have $|\bs|=0$.
For all $t \in [0,1] \setminus (\bs \cup \overline{\ns})$ we claim that
\begin{equation}
\label{convgrad2}
\nabla u_n(t) \to \nabla u(t) \quad \mbox{ strongly in }L^2(\Om;\R^2).
\end{equation}
In fact,
since $t \not \in \bigcup_{h=k}^\infty \bs_h$ for some $k$,
by Lemma \ref{convD} we have
\begin{equation*}
\lim_n \|\nabla u^{\delta_n}_{\varepsilon_n,a_n}(t)-\nabla u_{a_n}(t)\|=0;
\end{equation*}
for $t \not \in \overline{\ns}$, by Lemma \ref{extension} we have that
$\nabla u_{a_n}(t) \to \nabla u(t)$ strongly in $L^2(\Om; \R^2)$
and so \eqref{convgrad2} holds.
\par
Since $g_{\eps_n} \to g$ strongly in $W^{1,1}([0,1];H^1(\Om))$,
we deduce that for a.e. $\tau \in [0,1]$
\begin{equation*}
\nabla \dot{g}_{\eps_n}(\tau) \to \nabla \dot{g}(\tau)
\quad \mbox{ strongly in }L^2(\Om;\R^2).
\end{equation*}
Since $\Es_n(0) \to \Es(0)$ by \eqref{conv0} and
$o^{\delta_n}_{\eps_n} \to 0$, by semicontinuity of the energy and by
\eqref{discrenergybis} we have that for all
$t \in D$
\begin{equation}
\label{energyconv*}
\Es(t) \le \liminf_n \Es_n(t) \le \limsup_n \Es_n(t) \le
\Es(0)+2 \int_0^t \int_{\Om}
\nabla u(\tau) \nabla \dot{g}(\tau) \,dx\,d\tau.
\end{equation}
In view of \eqref{energybelow}, we conclude that for all $t \in D$
\begin{equation*}
\Es(t)=
\Es(0)+2 \int_0^t \int_{\Om}
\nabla u(\tau) \nabla \dot{g}(\tau) \,dx\,d\tau,
\end{equation*}
and since $\nabla u(\cdot)$ and $\hs^1(\Gamma(\cdot))$ are left continuous at $t \not \in D$ 
and so $\Es(\cdot)$ is,
we conclude that the equality holds for all $t \in [0,1]$. As a consequence
$\{t \to u(t)\,, t \in [0,1]\}$ is a quasi-static evolution of brittle fractures.
Let us prove that \eqref{energyconv*} is indeed true for all $t \in [0,1]$. In fact, if $t \notin D$,
it is sufficient to prove
\begin{equation}
\label{truelsc}
\liminf_n \Es_n(t) \ge \Es(t).
\end{equation}
Considering $s \ge t$ with $s \in D$, by \eqref{discrenergybists}
we have
$$
\Es_n(s) \le \Es_n(t)+
\int_{t^{\delta_n}_{j_n}}^{s^{\delta_n}_{j_n}} \int_{\Om} \nabla u_n(\tau)
\nabla \dot{g}_{\eps_n}(\tau)\,dx\,d\tau
+o^{\delta_n}_{\eps_n}
\quad
t^{\delta_n}_{j_n} \le t <t^{\delta_n}_{j_n+1},\;
s^{\delta_n}_{j_n} \le s <s^{\delta_n}_{j_n+1},
$$
so that
$$
\liminf_n \Es_n(t) \ge \Es(s)-
\int_t^s \int_{\Om} \nabla u(\tau) \nabla \dot{g}(\tau)\,dx\,d\tau.
$$
Letting $s \searrow t$, since $\Es(\cdot)$ is continuous, we have
\eqref{truelsc} holds. By \eqref{energyconv*} we deduce that $\Es_n(t) \to \Es(t)$ for all
$t \in [0,1]$, so that point $(b)$ is proved.
\par
We now come to point $(a)$.
Since $\lambda(t) \ge \hs^1(\Gamma(t))$ for all $t \in [0,1]$,
by \eqref{convgrad2} and point $(b)$, we deduce that $\lambda=\hs^1(\Gamma(\cdot))$ in $[0,1]$
up to a set of measure $0$. Since they are increasing functions,
we conclude that $\lambda$ and $\hs^1(\Gamma(\cdot))$
share the same set of continuity points $[0,1] \setminus \ns$, and that
$\lambda=\hs^1(\Gamma(\cdot))$ on $[0,1] \setminus \ns$. In view of \eqref{convgrad2},
point $(a)$ is thus established for all $t$ except $t \in (\bs \cup \overline{\ns}) \setminus \ns$.
In order to treat this case, we use the following argument.
Considering the measures $\mu_n:=\hs^1 \res \Gamma_n(t)$, we have
that, up to a subsequence, $\mu_n \weakst \mu$ weakly-star in the sense of measures,
and as a consequence of Ambrosio's Theorem we have $\hs^1 \res \Gamma(t) \le \mu$ as
measures. Since $t \notin \ns$ we have $\mu_n(\R^2) \to \hs^1(\Gamma(t))$, and so
we deduce
$\hs^1 \res \Gamma(t) =\mu$. Let us consider now $u_n(t)$; we have up to a subsequence
$u_n(t) \to u$ in $SBV(\Om)$ for some $u \in SBV(\Om)$. Setting $u_n(t):=g^{\delta_n}_{\eps_n}(t)$
and $u:=g(t)$ on $\Om_D$, we have $u_n(t) \to u$ in $SBV(\Om')$, and
as a consequence of Ambrosio's Theorem, we get that $\hs^1 \res \Sg{g(t)}{u}
\le \mu=\hs^1 \res \Gamma(t)$, that is
$\Sg{g(t)}{u} \subseteq \Gamma(t)$. By Theorem \ref{jumptransfer}, we deduce
that $u$ is a minimum for
$$
\min \{ \|\nabla v\|^2\,:\, \Sg{g(t)}{v} \subseteq \Gamma(t)
\mbox{ up to a set of $\hs^1$-measure $0$ } \},
$$
and by uniqueness of the gradient we get that $\nabla u=\nabla u(t)$, so that
the proof is concluded.
\end{proof}

\section{Piecewise Affine Transfer of Jump and Proof of Proposition \ref{pminpropD}}
\label{secmin}

The proof of Proposition \ref{pminpropD} is based on the following proposition,
which is a variant of Theorem \ref{jumptransfer} in the context of piecewise affine
approximation.

\begin{proposition}
\label{pctransfer}
Given $\eps_n \to 0$, let $g^r_n\in H^1(\Om)$ be such that $g^r_n \in \afenom$ and
$g^r_n \to g^r$ strongly in $H^1(\Om)$ for all $r=0,\dots,i$.
If $\unr \in \afenaom$ is such that $\unr \to \ur$ in $SBV(\Om)$ for $r=0,\dots,i$,
then for all $v \in SBV(\Om)$ with $\hs^1 \left( \Sg{g^i}{v} \right) <+\infty$ and $\nabla v \in
L^2(\Om;\R^2)$, there exists $v_n \in \afenaom$ such that $v_n \to v$ strongly in $L^1(\Om)$,
$\nabla v_n \to \nabla v$ strongly
in $L^2(\Om;\R^2)$ and
\begin{equation}
\label{pctransfjump}
\limsup_n \hs^1 \left( \Sg{g^i_n}{v_n} \setminus
\bigcup_{r=0}^i \Sg{g^r_n}{u^r_n}\right) \le
\mu(a) \hs^1 \left( \Sg{g(t)}{v} \setminus
\bigcup_{r=0}^i \Sg{g^r}{u^r} \right),
\end{equation}
where $\mu\,:\,]0;\frac{1}{2}[ \to \R$ with $\lim_{a \to 0^+} \mu(a)=1$.
\end{proposition}

In view of Proposition \ref{pctransfer}, we can now prove Proposition \ref{pminpropD}.

\begin{proof}[Proof of Proposition \ref{pminpropD}.]
Notice that, in order to prove \eqref{minpropD},
it is sufficient to prove the existence of $\mu\,:\,]0;\frac{1}{2}[ \to \R$ with
$\lim_{a \to 0^+} \mu(a)=1$ such that,
given
$t \in D$, for every $0=t_0 \le \ldots \le t_r \le \ldots \le t_i=t$, $t_r \in D$,
for all $v \in SBV(\Om)$ we have
\begin{equation}
\label{mineqk*}
\|\nabla u_a(t)\|^2 \le
\|\nabla v\|^2 +
\mu(a) \hs^1 \left( \Sg{g(t)}{v} \setminus
\bigcup_{r=0}^i \Sg{g(t_r)}{u_a(t_r)} \right).
\end{equation}
In fact, taking the sup over all
possible $t_0, \ldots, t_i$, we get \eqref{minpropD}.
\par
We apply Proposition \ref{pctransfer} considering $g^r_n:=g^{\delta_n}_{\eps_n}(t_r)$, $g^r:=g(t_r)$,
$u^r_n:=u^{\delta_n}_{\eps_n,a}(t_r)$, and $u^r:=u_a(t_r)$ for $r=0,\dots,i$. There exists
$\mu\,:\,]0;\frac{1}{2}[ \to \R$ with $\lim_{a \to 0^+} \mu(a)=1$ such that
for $v \in SBV(\Om)$, there exists $v_n \in \afenaom$ with $\nabla v_n \to \nabla v$ strongly
in $L^2(\Om;\R^2)$ and
\begin{multline*}
\limsup_n \hs^1 \left( \Sg{g^{\delta_n}_{\eps_n}(t)}{v_n} \setminus
\bigcup_{r=0}^i \Sg{g^{\delta_n}_{\eps_n}(t_r)}{u^{\delta_n}_{\eps_n,a}(t_r)} \right) \le \\
\le
\mu(a) \hs^1 \left( \Sg{g(t)}{v} \setminus
\bigcup_{r=0}^i \Sg{g(t_r)}{u_a(t_r)} \right),
\end{multline*}
Comparing $u^{\delta_n}_{\eps_n,a}(t)$ and $v_n$ by means of \eqref{piecemin2}, we obtain
\begin{multline}
\label{mminimality}
\|\nabla u^{\delta_n}_{\eps_n,a}(t) \|^2 \le \|\nabla v_n \|^2+
\hs^1 \left( \Sg{g^{\delta_n}_{\eps_n}(t)}{v_n} \setminus \Gamma^{\delta_n}_{\eps_n,a}(t) \right) \le \\
\le
\|\nabla v_n \|^2+
\hs^1 \left( \Sg{g^{\delta_n}_{\eps_n}(t)}{v_n} \setminus
\bigcup_{r=0}^i \Sg{g^{\delta_n}_{\eps_n}(t_r)}{u^{\delta_n}_{\eps_n,a}(t_r)} \right),
\end{multline}
so that, passing to the limit for $n \to +\infty$, we obtain that \eqref{mineqk*} holds. Moreover, we have that
choosing $v=u_a(t)$, and taking the limsup in \eqref{mminimality}, we get that $\nabla u^{\delta_n}_{\eps_n,a}(t)
\to \nabla u_a(t)$ strongly in $L^2(\Om;\R^2)$.
\end{proof}

The rest of the section is devoted to the proof of Proposition \ref{pctransfer}.
It will be convenient, as in Section \ref{convres}, to consider
$\Om_D$ polygonal open bounded subset of $\R^2$ such that
$\Om_D \cap \Om=\emptyset$ and $\partial \Om \cap \partial \Om_D= \partial_D \Om$ up to a finite number of
vertices; we set $\Om':= \Om \cup \Om_D \cup \partial_D \Om$. We suppose that $\rb_\eps$ can
be extended to a regular triangulation of $\Om'$ which we still indicate by
$\rb_\eps$.
\par
We need several preliminary results.
Let us set $\znr:=u^r_n-g^r_n$, and let
us extend $\znr$ to zero on $\Om_D$. Similarly, we set $\zr:=u^r-g^r$, and
we extend $\zr$ to zero on $\Om_D$.
\par
Let $\sigma>0$, and let
$C$ be the set of corners of $\partial_D \Om$.
Let us fix $G \subseteq \R$ countable and dense:
we recall that for all $r=0,\dots,i$ we have up to a set of $\hs^{1}$-measure zero
$$
S_{\zr}= \bigcup_{c_1,c_2 \in G} \partial^* E_{c_1}(r) \cap \partial^*E_{c_2}(r),
$$
where $E_c(r):=\{x \in \Om'\,:\, \zr(x) > c\}$ and $\partial^*$ denotes the
essential boundary (see \cite{AFP}). Let us consider
$$
J_j:=\{x \in \bigcup_{r=0}^i S_{\zr} \setminus C\,:\,
\ulp (x)-\ulm (x)>\frac{1}{j} \mbox{ for some } l=0, \ldots, i\},
$$
with $j$ so large that $\hs^1(\bigcup_{r=0}^i S_{\zr} \setminus J_j) \le \sigma$.
Let $U$ be a neighborhood of $\bigcup_{r=0}^i S_{\zr}$
such that $|U| \le \frac{\sigma}{j^2}$.
Following \cite[Theorem 2.1]{FL} (see Fig.3), we can find a finite
disjoint collection of closed cubes $\{Q_k\}_{k=1, \ldots,K}$
with center $x_k \in J_j$, edge of length $2r_k$ and
oriented as the normal $\nu(x_k)$ to $S_{\umk}$ at $x_k$,
such that $\bigcup_{k=1}^K Q_k \subseteq U$ and
$\hs^{1}(J_j \setminus \bigcup_{k=1}^K Q_k) \le \sigma$. Moreover
for all $k=1, \ldots,K$ there exists $r(k) \in \{0, \ldots, i\}$
and $\cksu,\cksd>0$ such that
$$
\hs^1 \left( \left[ \bigcup_{r=0}^i S_{\zr} \setminus S_{\umk}
\right] \cap Q_k \right) \le \sigma r_k,
$$
and the following hold

\begin{center}
\psfig{figure=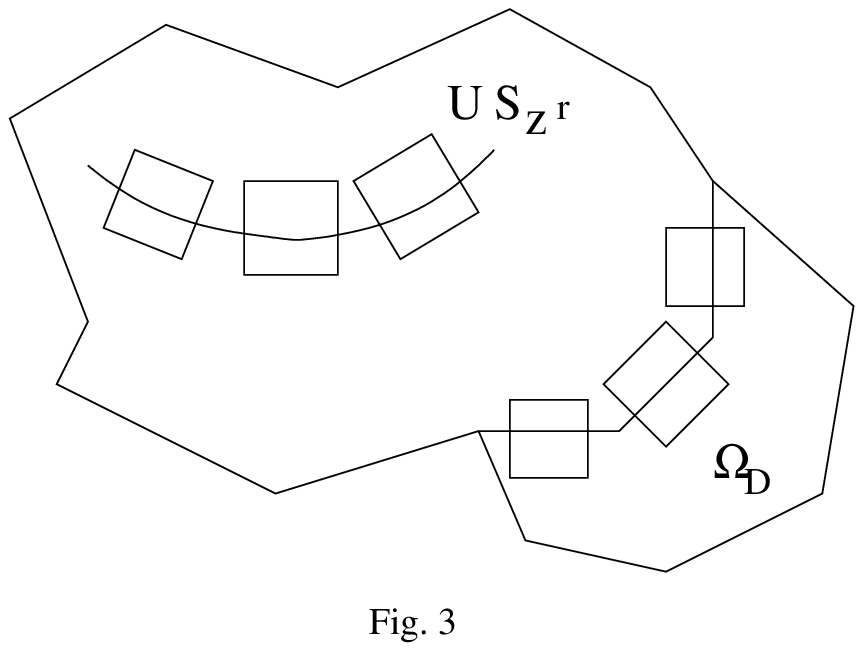}
\end{center}

\begin{itemize}
\item[(a)] if $x_k \in \Om$ then $Q_k \subseteq \Om$, and if $x_k \in \partial_D \Om$ then
$Q_k \cap \partial_D \Om=H_k$, where
$H_k$ denotes the intersection of $Q_k$ with the straight line
through $x_k$ orthogonal to $\nu(x_k)$;
\item[{}]
\item[(b)] $\hs^{1} (S_{\umk} \cap \partial Q_k)=0$;
\item[{}]
\item[(c)] $r_k \le c \hs^{1} (S_{\umk} \cap Q_k)$ for some $c>0$;
\item[{}]
\item[(d)] $\umkm(x) < \cksu < \cksd < \umkp(x)$ and
$\cksd-\cksu \ge \frac{1}{2j}$;
\item[{}]
\item[(e)] $\hs^{1} ([S_{\umk} \setminus
\partial^* E_{\cks}(r(k))] \cap Q_k)
\le \sigma r_k$ for $s=1,2$;
\item[{}]
\item[(f)] if $s=1,2$, $\hs^{1} (\{y \in \partial^* E_{\cks}(r(k)) \cap Q_k\,:\,
{\rm dist}(y,H_k) \ge \frac{\sigma}{2} r_k\}) <\sigma r_k$;
\item[{}]
\item[(g)] if $Q^+_k:= \{x\in Q_k \,|\, x \cdot \nu(x_k) >0\}$ and $s=1,2$
\begin{equation}
\label{density1}
\| 1_{E_{\cks}(r(k)) \cap Q_k} -1_{Q^+_k}\|_{L^1(\Om')}
\le \sigma^2 r_k^2;
\end{equation}
\item[{}]
\item[(h)] $\hs^{1} ((S_v \setminus S_{\umk}) \cap Q_k) < \sigma r_k$ and
$\hs^{1} (S_v \cap \partial Q_k)=0$.
\end{itemize}

Let us indicate by $R_k$ the intersection of $Q_k$
with the strip centered in $H_k$ with width $2\sigma r_k$, and let
us set $V_k^\pm:=\{x_k\pm r_k e(x_k)+s \nu(x_k):\, s \in \R\} \cap R_k$, where
$e(x_k)$ is such that $\{e(x_k), \nu(x_k)\}$ is an orthonormal base of $\R^2$ with the
same orientation of the canonical one.
\par
For all $B \subseteq \Om'$, let us set
\begin{equation*}
\rs_n(B):=\{T \in \rb_{\eps_n} \,:\, T \cap B \not= \emptyset\}, \quad \quad
\ts^k_n(B):=\{T \in \tb(\umkn)\,:\, T \cap B \not= \emptyset\}.
\end{equation*}
In the following, we will often indicate with the same symbol a family of
triangles and their support in $\R^2$, being clear from the context
in which sense has to be intended.
We will consider $\umkn$ defined pointwise in
$\Om' \setminus \overline{S}_{\umkn}$ and so the upper levels of $\umkn$ are
intended as subsets of $\Om' \setminus \overline{S}_{\umkn}$.

\begin{lemma}
\label{upperlevels}
For all $k=1, \ldots,K$ there exists
$c_n^k \in [\cksu,\cksd]$
such that, setting $\Enk:=\{x \in \rskn\,:\,
\umkn(x) >c^k_n\}$, we have
\begin{equation}
\label{coarea1}
\limsup_n \sum_{k=1}^K
\hs^1 \left( \big( \partial_{\rskn} \Enk
\big) \setminus S_{\umkn} \right)=o_\sigma,
\end{equation}
and
\begin{equation}
\label{nearset1}
\limsup_n \|1_{\Enk}-1_{Q_k^+}\|_{L^1(\Om')} \le
\sigma^2 r_k^2,
\end{equation}
where $\partial_{\rskn}$ denotes the boundary operator in
$\rskn$, and $o_\sigma \to 0$ as $\sigma \to 0$.
\end{lemma}

\begin{proof}
Note that for $n$ large we have $\bigcup_{k=1}^K \rskn \subseteq U$, so that
$|\bigcup_{k=1}^K \rskn| \le \frac{\sigma}{j^2}$. By H\"older inequality and since
$\|\nabla \znr \| \le C'$ for all $r=0,\dots,i$, it follows that
\begin{equation*}
\sum_{r=0}^i \int_{\{\cup_k \rskn: r(k)=r\}} |\nabla \znr| \,dx \le
\sum_{r=0}^i \|\nabla  \znr\| \frac{\sqrt{\sigma}}{j} \le (i+1)C'\frac{\sqrt{\sigma}}{j}.
\end{equation*}
Following \cite[Theorem 2.1]{FL}, we can apply coarea-formula for BV-functions (see \cite{AFP})
taking into account that $\umkn$ belongs to $SBV(\Om')$ so that the singular
part of the derivative is carried only by $S_{\umkn}$:
since for $n$ large the $\rskn$'s are disjoint, we obtain
\begin{equation}
\label{coarea}
\sum_{k=1}^K \int_\R
\hs^1 \Big( \big( \partial E_{c,n}(r(k))
\cap \rskn \big) \setminus S_{\umkn} \Big)\, dc \le (i+1)C' \frac{\sqrt{\sigma}}{j},
\end{equation}
where $E_{c,n}(r(k)):=\{x \in \Om' \setminus \overline{S}_{\umkn}\,:\,
\umkn(x) >c\}$, and so
\begin{equation*}
\sum_{k=1}^K \int_{\cksu}^{\cksd}
\hs^1 \Big( \big( \partial E_{c,n}(r(k))
\cap \rskn \big) \setminus S_{\umkn} \Big) \, dc
\le (i+1)C' \frac{\sqrt{\sigma}}{j}.
\end{equation*}
Notice that we can use the topological boundary instead of the reduced boundary of $E_{c,n}(r(k))$
in \eqref{coarea} since $\umkn$ is piecewise affine, and so $\partial E_{c,n}(r(k)) \setminus
\partial^* E_{c,n}(r(k)) \not= \emptyset$ just for a finite number of $c$'s.
By the Mean Value  Theorem we have that there exist
$c_n^k \in [\cksu,\cksd]$ such that
\begin{equation*}
\sum_{k=1}^K
\hs^1 \Big( \big( \partial E_{c_n^k,n}(r(k))
\cap \rskn \big) \setminus S_{\umkn} \Big) \, \le 2iC' \sqrt{\sigma},
\end{equation*}
and taking the limsup for $n \to +\infty$, we get \eqref{coarea1}.
Let us come to \eqref{nearset1}. Since
$$E_{\cksd,n}(r(k)) \subseteq E_{c_n^k,n}(r(k)) \subseteq E_{\cksu,n}(r(k)),$$ by
\eqref{density1} we have that for $n$ large
$$
\|1_{E_{c_n^k,n}(r(k)) \cap Q_k} -1_{Q_k^+}\|_{L^1(\Om')} \le \sigma^2 r_k^2,
$$
and so, since $|\rskn \setminus Q_k| \to 0$, we conclude that
\eqref{nearset1} holds.
\end{proof}

Fix $k \in \{1, \ldots, K\}$, and let us consider the family
$\ts^k_n(\Enk)$. Let us modify this family in the following way.
Let $T \in \ts^k_n(\Enk)$; we keep it if
$|T \cap \Enk|> \frac{1}{2}|T|$, and we erase it
otherwise. Let $\Enkp$ be this new family of triangles, and let
$\Enkm$ be its complement in $\ts^k_n(\rskn)$.

\begin{lemma}
\label{en+}
For all $k=1, \ldots,K$ we have
\begin{equation}
\label{estper}
\limsup_n \sum_{k=1}^K
\hs^1 \left( \partial_{\rskn} \Enkp \setminus S_{\umkn} \right)
=o_\sigma,
\end{equation}
and
\begin{equation}
\label{estarea}
\limsup_n \|1_{\Enkp}-1_{Q_k^+}\|_1 \le 4\sigma^2 r_k^2,
\end{equation}
where $o_\sigma \to 0$ as $\sigma \to 0$.
\end{lemma}

\begin{proof}
Let $T \in \ts^k_n(\Enk)$. Since $\umkn$ is affine on
$T$, it follows that $T \cap \Enk$ is either a triangle with at least two edges
contained in the edges of $T$ or a trapezoid with three edges contained
in the edges of $T$. Let $l(T)$ be the edge inside $T$ where $\umkn=c_n^k$,
where $c^k_n$ is the value determining $\Enk$ (we consider $l(T)=\emptyset$ if
${\rm int}(T) \subseteq \Enk$).
In the case $T \in \Enkp$ as in the case $T \in \Enkm$, since the angles
of the triangles of $\tb(\umkn)$ are uniformly bounded away from $0$ and from
$\pi$, arguing as in Lemma \ref{simplecurve}, we deduce that keeping or erasing $T$,
we increase $\partial_{\rskn} \Enk$ of a quantity which is less than
$c\hs^1(l(T))$ with $c$ independent of $\varepsilon_n$. Then we have
$$
\sum_{k=1}^K \hs^1(\partial_{\rskn} \Enkp \setminus \partial_{\rskn} \Enk)
\le \sum_{k=1}^K
\sum_{T \in \ts^k_n(\Enk)} c\hs^1(l(T))
\le c \sum_{k=1}^K \hs^1( \partial_{\rskn} \Enk \setminus S_{\umkn}),
$$
so that taking the limsup for $n \to +\infty$ and in view of \eqref{coarea1}
we deduce that \eqref{estper} holds.
\par
Let us come to \eqref{estarea}. Note that $|\ts^k_n(\partial Q_k^+)| \to 0$
as $n \to +\infty$. Then if
$A^{k,+}_n:=\{T \in \tb(\umkn)\,:\, T \subseteq \tint{Q_k^+}\}$, for $n$
large we have
$$
|Q_k^+ \setminus \Enkp| \le |A^{k,+}_n \setminus \Enkp|+
|\ts^k_n(\partial Q_k^+)| \le 2|Q_k^+ \setminus \Enk|+
|\ts^k_n(\partial Q_k^+)|,
$$
where the last inequality follows by construction of $\Enkp$.
Taking the limsup for $n \to +\infty$, in view of \eqref{nearset1} we get
$$
\limsup_n |Q_k^+ \setminus \Enkp| \le 2\sigma^2 r_k^2.
$$
The inequality $\limsup_n |\Enkp \setminus Q_k^+| \le 2\sigma^2 r_k^2$ follows
analogously.
\end{proof}

For all $k=1, \ldots, K$ and $s \in \R$, let us set
$$
\hkn(s) := \{x + s \nu(x_k),\, x\in H_k\}.
$$

\begin{lemma}
\label{fubini}
There exist $s^+_n \in ]\frac{\sigma}{4}r_k, \frac{\sigma}{2} r_k[$ and
$s^-_n \in ]-\frac{\sigma}{2} r_k, -\frac{\sigma}{4}r_k[$ such that, setting
$\hknp:=\hkn(s^+_n)$ and $\hknm:=\hkn(s^-_n)$ we have for $n$ large enough
\begin{equation*}
\hs^1(\hknp \setminus \Enkp) \le 20\sigma r_k,
\quad
\hs^1(\hknm \cap \Enkp) \le 20\sigma r_k.
\end{equation*}
\end{lemma}

\begin{proof}
By (\ref{estarea}) we can write for $n$ large
$$
\int_{\frac{\sigma}{4}r_k}^{\frac{\sigma}{2} r_k}
\hs^1(H_k(s) \setminus \Enkp) \,ds \le
5 \sigma^2 r_k^2,
$$
so that we get $s^+_n \in ]\frac{\sigma}{4}r_k, \frac{\sigma}{2} r_k[$
with
$$
\hs^1(H_k(s^+_n) \setminus \Enkp) \le 20\sigma r_k.
$$
Similarly we can reason for $s^-_n$.
\end{proof}

Let $S^{k,+}_n$ be the straight line containing
$\hknp$: up to replacing $\hknp$ by
the connected component of $S^{k,+}_n \cap \rskn$
to which it belongs, we may suppose that
$\hknp \setminus \Enkp$ is a finite union of segments
$l_j^+$ with extremes $A_j$ and $B_j$ belonging to
the edges of the triangles of $\ts_n^k(\rskn)$ such that for $n$ large
$$
\hs^1(\hknp \setminus \Enkp)=\hs^1 \left( \bigcup_{j=1}^{m} l_j^+ \right) \le 20\sigma r_k.
$$
By Lemma \ref{simplecurve},
for all $j$ there exists a curve $L_j^+$
inside the edges of the triangles of $\ts_n^k(\rskn)$ joining $A_j$ and $B_j$
and such that
\begin{equation}
\label{estsimplecurve}
\hs^1(L_j^+) \le c\hs^1(l_j^+),
\end{equation}
with $c$ independent of $\varepsilon_n$.
Let us set
$$
\gamma^{k,+}_n:= L_1^+ \cup B_1 A_2 \cup L_2^+ \cup \cdots \cup B_{m-1}A_m \cup
L_{m}^+.
$$
Similarly, let us construct $\gamma^{k,-}_n$ relative to $\hknm \cap E_n^+$.
Note that for $n$ large enough $\gamma^{k,+}_n \cap H_k(\sigma)=\emptyset$,
$\gamma^{k,-}_n \cap H_k(-\sigma)=\emptyset$, and
$\gamma^{k,+}_n \cap \gamma^{k,-}_n=\emptyset$.
Let us consider the connected component
$\cs_k^+$ of $\rskn \setminus \gamma^{k,+}_n$
containing $H_k(\sigma)$.
Similarly, let us consider the connected component $\cs_k^-$ of
$\rskn \setminus \gamma^{k,-}_n$ containing $H_k(-\sigma)$.
For $n$ large enough, by \eqref{estsimplecurve}
\begin{equation}
\label{estsimplecurve2}
\hs^1 \left( \partial_{\rskn} \cs_k^+
\setminus
\bigcup_{i=1}^{m-1} B_i A_{i+1}\right)
\le c \sum_{j=1}^{m} \hs^1(l_j^+) \le 20c\sigma r_k.
\end{equation}
A similar estimate holds for $\partial_{\rskn} \cs_k^-$.
\par
Let $\enkpt$ be the family of triangles obtained adding to $\Enkp$
those $T \in \Enkm$ such that $T \subseteq \cs_k^+$, and subtracting
those $T \in \Enkp$ such that $T \subseteq \cs_k^-$. Let $\enkmt$
be the complement of $\enkpt$ in $\ts_n^k(\rskn)$.
\par
We claim that there exists
$C>0$ independent of $n$ such that for all $k=1, \ldots, K$ and for $n$ large
\begin{equation}
\label{estper3}
\hs^1 \left( \partial_{\rskn} \enkpt \setminus  \partial_{\rskn} \Enkp \right) \le
C\sigma r_k.
\end{equation}
In fact, let $\zeta$ be an edge of $\partial_{\rskn} \enkpt \setminus
\partial_{\rskn} \Enkp$, that is $\zeta$ belongs to a triangle $T$ that has
been changed in the operation above described.
Let us assume for instance that $T \in \Enkm$ and $T \subseteq \cs_k^+$.
If $T'$ is such that $T \cap T'=\zeta$, then $T' \in \Enkm$: in fact
if by contradiction $T' \in \Enkp$, then $T' \in \enkpt$ and so
we would have $\zeta \not\in \partial_{\rskn} \enkpt$ which is absurd.
Similarly we get $T' \not \subseteq \cs_k^+$. This means that
$\zeta \subseteq \partial_{\rskn} \cs_k^+$, and
since the horizontal edges of $\gamma^{k,+}_n$ intersect by construction only elements of $\Enkp$,
we deduce that $\zeta \subseteq \partial_{\rskn} \cs_k^+ \setminus
\left( \cup_{i=1}^{m} A_iB_i \right)$,
and by \eqref{estsimplecurve2} we conclude that \eqref{estper3} holds.
\par
We can summarize the previous results as follows.

\begin{lemma}
\label{families}
For all $k=1, \ldots, K$ there exist two families $\enkpt$ and $\enkmt$
of triangles with $\ts_n^k(\rskn)=\enkpt \cup \enkmt$,
$Q_k^+ \setminus R_k \subseteq \enkpt$ and $Q_k^- \setminus R_k \subseteq \enkmt$,
and such that
\begin{equation}
\label{estper2}
\limsup_n \sum_{k=1}^K \hs^1 \left( \partial_{\rskn} \enkpt \setminus S_{\umkn} \right)=
o_\sigma,
\end{equation}
where $o_\sigma \to 0$ as $\sigma \to 0$.
Moreover, in the case $x_k \in \partial_D \Om$, we can modify
$\enkpt$ or $\enkmt$ in such a way that $\enkpt \subseteq \Om$ or $\enkmt \subseteq \Om$.
\end{lemma}

\begin{proof}
We have that \eqref{estper2} follows from \eqref{estper} and
\eqref{estper3}, and the fact that $\sum_{k=1}^K r_k \le c$, with
$c$ independent of $\sigma$.
Let us consider the case $x_k \in \partial_D \Om$
with $Q_k^+ \setminus R_k \subseteq \Om$ (the other case being similar).
From \eqref{estper2} we have that for $n$ large
$\sum_{k=1}^K \hs^1 \left( \partial_{\rskn} \enkpt \cap Q_k^- \right) \le o_\sigma$
because $\umkn=R_{\eps_n} g_{h_n}(r(k))$ on $Q_k^-$ and so there are no jumps
in $Q_k^-$. We can thus redefine $\enkpt$ subtracting those triangles that are
in $Q_k^-$ obtaining again \eqref{estper2}.
\end{proof}

We are now in position to prove Proposition \ref{pctransfer}.

\begin{proof}[Proof of Proposition \ref{pctransfer}.]
We work in the context of $\Om'$. For all $v \in SBV(\Om')$ with $v=g^i$ on $\Om_D$,
$\hs^1(S_v) <+\infty$ and $\nabla v \in L^2(\Om';\R^2)$,
we have to construct $v_n \in SBV(\Om')$
such that $v_n=g^i_n$ on $\Om_D$, $(v_n)_{|\Om} \in \afenaom$, $v_n \to v$ strongly in $L^1(\Om')$,
$\nabla v_n \to \nabla v$ strongly in
$L^2(\Om';\R^2)$ and
\begin{equation}
\label{mineqk}
\limsup_n \hs^1 \left( S_{v_n} \setminus
\bigcup_{r=0}^i S_{u^r_n}\right) \le
\mu(a) \hs^1 \left( S_{v} \setminus
\bigcup_{r=0}^i S_{u^r} \right),
\end{equation}
where we suppose that
$u^r_n$ and $u^r$ are extended to $\Om'$ setting $u^r_n:=g^r_n$, and $u^r:=g^r$ on $\Om_D$ respectively.
\par
We set $v=g^i+w$, where $w \in SBV(\Om')$ with $w=0$ on $\Om_D$. By density, it is sufficient
to consider the case $w \in L^\infty(\Om')$.
Up to reducing $U$, we may assume that $\|\nabla g^i\|_{L^2(U;\R^2)} < \sigma$
and $\|\nabla w\|_{L^2(U;\R^2)} <\sigma$. Let $R'_k$ be a rectangle centered in $x_k$,
oriented as $R_k$, and such that $\overline{R'_k} \subset {\rm int}R_k$ and
$\hs^1(S_w \cap (R_k \setminus R'_k)) < \sigma r_k$.
We claim that there exists $\wsig \in SBV(\Om')$ with $\wsig=w$ on $\bigcup_{k=1}^K R'_k$ and
$\wsig=0$ in $\Om_D$ such that
\begin{itemize}
\item[]
\item[(1)] $\|w -\wsig\|+\|\nabla w-\nabla \wsig\| \le \sigma$;
\item[]
\item[(2)] $\hs^1(S_{\wsig} \cap (Q_k \setminus R'_k)) \le o_\sigma r_k$, with $o_\sigma \to 0$ as
$\sigma \to 0$;
\item[]
\item[(3)] $\hs^1(S_{\wsig} \setminus \bigcup_{k=1}^K R_k) \le \hs^1(S_w \setminus \bigcup_{k=1}^K R_k)+\sigma$;
\item[]
\item[(4)] $S_{\wsig} \setminus \bigcup_{k=1}^K R_k$ is union of disjoint segments
with closure contained in $\Om \setminus \bigcup_{k=1}^K R_k$;
\item[]
\item[(5)] $\wsig$ is of class $W^{2,\infty}$ on $\Om \setminus \left( \bigcup_{k=1}^K R_k \cup
\overline{S_{\wsig}} \right)$.
\end{itemize}
In fact, by Proposition \ref{regularization}, there exists
$w_m \in SBV(\Om')$ with $w_m=0$ in $\Om' \setminus \Omb$
such that $w_m \to w$
strongly in $L^2(\Om')$, $\nabla w_m \to \nabla w$ strongly in
$L^2(\Om';\R^2)$, $S_{w_m}$ is polyhedral with
$\overline{S_{w_m}} \subseteq \Om$,
$w_m$ is of class $W^{2,\infty}$ on $\Om \setminus \left( \bigcup_{k=1}^K R_k \cup
\overline{S_{w_m}} \right)$, and
$\lim_m \hs^{1}(A \cap S_{w_m})=\hs^{1}(A \cap S_w)$
for all $A$ open subset of $\Om'$ with $\hs^{1}(\partial A \cap S_{w})=0$.
It is not restrictive to assume that $\hs^1(S_{w} \cap \partial R_k)=0$ and
$\hs^1(S_{w_m} \cap \partial R_k)=0$ for all $m$.
Let $\psi_k$ be a smooth function such that $0 \le \psi_k \le 1$, $\psi_k=1$ on
$R'_k$ and $\psi_k=0$ outside $R_k$. Setting $\psi:=\sum_{k=1}^K \psi_k$, let
us consider $\tilde{w}_m:= \psi w+ (1-\psi) w_m$. Note that
$\tilde{w}_m \to w$ strongly in $L^2(\Om')$, $\nabla \tilde{w}_m \to \nabla w$ strongly in
$L^2(\Om';\R^2)$, $\tilde{w}_m=0$ in $\Om_D$. Moreover,
by capacity arguments, we may assume that $S_{\tilde{w}_m} \setminus \bigcup_{k=1}^K R_k$
is a finite union of disjoint segments with closure contained in
$\Om \setminus \bigcup_{k=1}^K R_k$.
Finally, for $m \to +\infty$, we have
$$
\hs^1(S_{\tilde{w}_m} \setminus \bigcup_{k=1}^K R_k) \to
\hs^1(S_{w} \setminus \bigcup_{k=1}^K R_k)),
$$
$$
\hs^1(S_{\tilde{w}_m} \cap \bigcup_{k=1}^K (Q_k \setminus R_k)) \to
\hs^1(S_{w} \cap \bigcup_{k=1}^K (Q_k \setminus R_k))
$$
and
$
\limsup_m \hs^1(S_{\tilde{w}_m} \cap (R_k \setminus R_k')) \le
2\hs^1(S_w \cap (R_k \setminus R_k')) \le 2\sigma r_k.
$
Then we can take $\wsig:=\tilde{w}_m$ for $m$ large enough.
\par
Let $S_{\wsig} \setminus \bigcup_{k=1}^K Q_k:= \bigcup_{j=1}^m l_j$, where, by capacity arguments,
we can always assume that $l_j$ are
disjoint segments with closure contained in $\Om \setminus \bigcup_{k=1}^K Q_k$. We define a triangulation
$\tb_n \in \ts_{\varepsilon_n,a}(\Om')$ specifying its adaptive vertices
as follows.
Let us consider the families $\rskn$ and $\rs_n(l_j)$ for $k=1, \ldots, K$
and $j=1, \ldots, m$. Note that for $n$ large enough,
$\rs_n(Q_{k_1}) \cap \rs_n(Q_{k_2})= \emptyset$ for
$k_1 \not= k_2$, $\rs_n(l_{j_1}) \cap \rs_n(l_{j_2})= \emptyset$ for $j_1 \not= j_2$,
and $\rs_n(Q_k) \cap \rs_n(l_j)= \emptyset$ for every $k,j$.
We consider inside the triangles of $\rs_n(Q_k)$ the adaptive vertices
of $\tb(\umkn)$.
Passing to $\rs_n(l_j)$, by density arguments it is not restrictive to assume that $l_j$
does not pass through the vertices of $\rb_{\eps_n}$ and that its extremes belong to the
edges of $\rb_{\eps_n}$.
Let $\zeta:=[x,y]$ be an edge of $\rs_n(l_j)$ such that $l_j \cap \zeta =\{P\}$.
Proceeding as in \cite{N}, we take as adaptive vertex of $\zeta$
the projection of $P$ on $\{tx+(1-t)y:\, t \in [a,(1-a)]\}$.
Connecting these adaptive vertices, we obtain
an {\it interpolating} polyhedral curve $\tilde{l}_j$ with
\begin{equation}
\label{estinterp}
\hs^1(\tilde{l}_j) \le \mu(a) \hs^1(l_j),
\end{equation}
where $\mu$ is an increasing function such that $\lim_{a \to 0} \mu(a)=1$.
Finally, in the remaining edges, we can consider any admissible adaptive vertex,
for example the middle point.
\par
Let us define $w_n \in SBV(\Om')$ in the following way. For all
$Q_k$, let $w_n$ be equal to $\wsig$ on $\rs_n(Q_k) \setminus R_k$,
equal to the reflection of ${\wsig}_{|Q_k^+ \setminus R_k}$ with respect to $H_k(\sigma)$ on $\enkpt \cap R_k$ and
equal to the reflection of ${\wsig}_{|Q_k^- \setminus R_k}$ with respect to $H_k(-\sigma)$ on $\enkmt \cap R_k$, where
$\tilde{E}^{k,\pm}_n$ are defined as in Lemma \ref{families}.
On the other elements of $\tb_n$, let us set $w_n=\wsig$.
Notice that $w_n=0$ on $\Om_D$ and that inside each $\rskn$, all the discontinuities of
$w_n$ are contained in $\partial_{\rskn} \enkpt \cup V_k \cup P_{\wsig}^k$, where $P_{\wsig}^k$ is
the union of the polyhedral jumps of $\wsig$ in $\rs_n(Q_k)$ and of their reflected version
with respect to $H_k(\pm \sigma)$.
By Lemma \ref{families} and since $\sum_{k=1}^K \hs^1(V_k \cup P_{\wsig}^k) \le o_\sigma$ with
$o_\sigma \to 0$ as $\sigma \to 0$,
and $\hs^1 \left( \bigcup_{r=0}^i S_{\zr} \setminus \bigcup_{k=1}^K Q_k \right)
\le 2\sigma$, we have that
\begin{multline*}
\limsup_n \hs^1 \left( S_{w_n} \setminus \bigcup_{r=0}^i S_{\znr} \right) \le \\
\le \hs^1 \left( S_{\wsig} \setminus \bigcup_{k=1}^K Q_k \right) +
\limsup_n \hs^1 \left( (S_{w_n} \setminus \bigcup_{r=0}^i S_{\znr}) \cap \rskn \right) \le \\
\le \hs^1 \left( S_{\wsig} \setminus \bigcup_{r=0}^i S_{\zr} \right)+
\hs^1 \left( \bigcup_{r=0}^i S_{\zr} \setminus \bigcup_{k=1}^K Q_k \right)+
o_\sigma \le \\
\le \hs^1 \left( S_{\wsig} \setminus \bigcup_{r=0}^i S_{\zr} \right)+o_\sigma,
\end{multline*}
and since $\|\nabla \wsig\|_{L^2(U;\R^2)} \le o_\sigma$ we get for $n$ large
\begin{equation}
\label{estimategrad}
\|\nabla w_n\|^2_{L^2(\bigcup_{k=1}^K \rskn)} \le o_\sigma.
\end{equation}
We now want to define an interpolation $\twn$ of $w_n$ on $\tb_n$.
Firstly, we set $\twn=0$ on all regular triangles of $\Om_D$.
Passing to the triangles in $\rskn$ (see fig.4), by Lemma \ref{guscio},
we know that for $n$ large enough, we have
$$
\hs^1(\partial \rs_n(V_k)) \le c \hs^1(V_k), \quad
\hs^1(\partial \rs_n(P_{\wsig}^k)) \le c \hs^1(P_{\wsig}^k),
$$
with $c$ independent of $n$.
If $T \in \rs_n(V_k) \cup \rs_n(P_{\wsig}^k)$, we set $\twn=0$ on $T$;
otherwise, we define $\twn$ on $T$ as the affine interpolation of $w_n$.
\par
Since $\nabla \twn$ is uniformly bounded on $\rs_n(H_k(\pm \sigma))$,
$|\rs_n(H_k(\pm \sigma))| \to 0$ and since $w_n$ is uniformly bounded in $W^{2,\infty}$ on
the triangles contained in $\rs_n(Q_k) \setminus \rs_n(V_k \cup P_{\wsig}^k \cup H_k(\pm \sigma))$
we have by the interpolation estimate \eqref{interp1} and by \eqref{estimategrad}
\begin{equation}
\label{estgradqint}
\limsup_n \|\nabla \twn\|^2_{L^2(\bigcup_{k=1}^K \rskn)} \le o_\sigma.
\end{equation}
Moreover we have
\begin{equation}
\label{estjumpqint}
\limsup_n \sum_{k=1}^K
\hs^1 \left( \Big( S_{\twn} \setminus \bigcup_{r=0}^k S_{\znr} \Big) \cap \rskn \right) \le o_\sigma.
\end{equation}
Let us come to the triangles not belonging to $\rskn$ for $k=1, \ldots, K$.
For all $j=1, \ldots, m$, we denote by $\hat{\rs}_n(l_j)$ the family of regular
triangles that have edges in common with triangles of $\rs_n(l_j)$.
For $n$ large we have that $\hat{\rs}_n(l_{j_1}) \cap \hat{\rs}_n(l_{j_2})= \emptyset$ for
$j_1 \not= j_2$. On every regular triangle $T \not \in \bigcup_{k=1}^K \rskn \cup
\bigcup_{j=1}^m \hat{\rs}_n(l_j)$,
we define $\twn$ as the affine interpolation of $\wsig$.
Since $\wsig$ is of class $W^{2,\infty}$ on $T$ and $T$ is regular,
we obtain by the interpolation estimate \eqref{interp1}
\begin{equation}
\label{estgradinside}
\|\twn - \wsig\|^2_{W^{1,2}(T)} \le \cint \varepsilon_n \|\wsig\|_{W^{2,\infty}}.
\end{equation}
Let us consider now those triangles that are contained in the elements
of $\bigcup_{j=1}^m \hat{\rs}_n(l_j)$. Following \cite{N}, we can define $\twn$
on every $T$ in such a way that $\twn$ admits discontinuities
only on $\tilde{l}_j$, and $\|\nabla \twn\|_{L^\infty(T)} \le \|\nabla \wsig\|_\infty$.
Since $|\hat{\rs}_n(l_j)| \to 0$ as $n \to \infty$, we deduce that
\begin{equation}
\label{estgradjumps}
\lim_n \|\nabla \twn\|^2_{L^2(\hat{\rs}_n(l_j))}=0.
\end{equation}
Moreover by \eqref{estinterp} and since
$\hs^1 \left( \bigcup_{r=0}^i S_{\zr} \setminus \bigcup_{k=1}^K Q_k \right) \le 2\sigma$,
we have
\begin{equation}
\label{estjumpinside}
\hs^1 \left( S_{\twn} \cap \bigcup_{j=1}^m \hat{\rs}_n(l_j) \right) \le \mu(a)
\hs^1 \left( S_{\wsig} \setminus \bigcup_{i=1}^k S_{\zr} \right) +o_\sigma,
\end{equation}
where $o_\sigma \to 0$ as $\sigma \to 0$.

\begin{center}
\psfig{figure=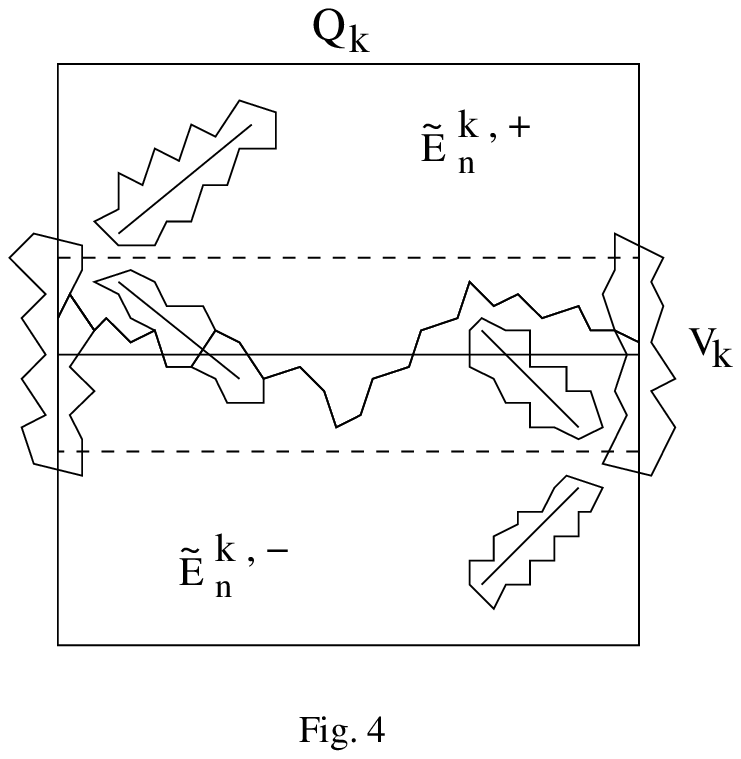}
\end{center}

We are now ready to conclude. Let us consider $\hwn \in \afenaom$ defined
as $\hwn:=g^i_n+\twn$. We have $\hwn \to g^i+\wsig$ strongly in $L^2(\Om')$.
By \eqref{estgradqint}, \eqref{estgradinside}, \eqref{estgradjumps} we get
$$
\limsup_n \|\nabla \hwn\|^2 \le \|\nabla g^i+\nabla \wsig\|^2 +o_\sigma,
$$
while by \eqref{estjumpqint} and \eqref{estjumpinside} we have
$$
\limsup_n \hs^1 \left( S_{\hwn} \setminus \bigcup_{r=0}^i S_{\znr} \right) \le \mu(a)
\hs^1 \left( S_{\wsig} \setminus \bigcup_{r=0}^i S_{\zr} \right)+o_\sigma.
$$
Letting now $\sigma \to 0$, using a diagonal argument, we conclude that Proposition \ref{pctransfer} holds.
\end{proof}

\section{Revisiting the approximation by Francfort and Larsen}
\label{remark}

In this section we show how the arguments of Section \ref{convres} may be used
to deal with the discrete in time approximation of quasi-static growth of
brittle fractures proposed by Francfort and Larsen in \cite{FL}.
More precisely, we prove that there is strong convergence of the gradient of the
displacement
(in particular convergence of the bulk energy) and convergence of
the surface energy at all times of continuity of the length of the crack;
moreover there is convergence of the total energy at any time.
\par
We briefly recall the notation employed in \cite{FL}. Let $I_\infty$ be countable and
dense in $[0,1]$, and let $I_n:=\{0=t^n_0 \le \ldots \le t^n_n=1\}$ such that
$(I_n)$ is an increasing sequence of sets whose union is $I_\infty$. Let $\Om \subseteq
\R^N$ be a Lipschitz bounded domain, and let $\partial \Om=\partial \Om^c_f \cup \partial \Om_f$,
where $\partial \Om^c_f$ is open in the relative topology. Let $\Om' \subseteq \R^N$ be open and
such that $\overline{\Om} \subseteq \Om'$, and let $g \in W^{1,1}([0,1];H^1(\Om'))$. At any
time $t^n_k$, Francfort and Larsen consider $u^n_k$ minimizer of
$$
\int_{\Om} |\nabla v|^2 \,dx+ \hs^{N-1} \left( S_v \setminus \left[
\bigcup_{0 \le j \le k-1} S_{u^n_j} \cup \partial \Om_f \right] \right)
$$
in $\{v \in SBV(\Om'):v=g(t^n_k)\mbox{ in }\Om'\setminus \overline{\Om}\}$.
Setting $u^n(t):=u^n_k$ for $t \in [t^n_k,t^n_{k+1}[$, and $\Gamma^n(t):=
\bigcup_{s \le t, s \in I_n} S_{u^n(s)} \cup \partial \Om_f$, they prove that
\begin{equation}
\label{estabovec}
\Es^n(t) \le \Es^n(0)+2\int_0^{t^n_k} \int_{\Om} \nabla u^n(\tau) \nabla \dot{g}(\tau) \,dx\,d\tau
+o_n, \quad\quad t\in [t^n_k,t^n_{k+1}[,
\end{equation}
where $\Es^n(t):=\int_{\Om} |\nabla u^n(t)|^2 \,dx+ \hs^{N-1}\left( \Gamma^n(t) \right)$ and
$o_n \to 0$ as $n \to +\infty$.
Using Theorem \ref{jumptransfer}, they obtain a subsequence of $(u^n(\cdot))$, still denoted by the same symbol,
such that $u^n(t) \to u(t)$ in $SBV(\Om')$ and $\nabla u^n(t) \to \nabla u(t)$ strongly in $L^2(\Om';\R^N)$
for all $t \in I_\infty$, with $u(t)$ a minimizer of
$$
\int_{\Om} |\nabla v|^2 \,dx+ \hs^{N-1} \left( S_v \setminus \Gamma(t) \right),
$$
where $\Gamma(t):= \bigcup_{s \in I_{\infty}, s \le t} S_{u(s)} \cup \partial \Om_f$.
The evolution $\{t \to u(t),\,t \in I_\infty\}$ is extended to the whole $[0,1]$ using the
approximation from the left in time.
\par
We can now use the arguments of Section \ref{convres}. Following Lemma \ref{extension}, 
it turns out that for all $t \in [0,1]$
\begin{equation}
\label{estbelowc}
\Es(t) \ge \Es(0)+2\int_0^t \int_\Om \nabla u(\tau) \nabla \dot{g}(\tau) \,dx\,d\tau.
\end{equation}
Moreover, by the Transfer of Jump and the uniqueness argument of Lemma \ref{aextension},
we have that $\nabla u^n(t) \to \nabla u(t)$ strongly
in $L^2(\Om';\R^N)$ for all $t \not\in \ns$, where $\ns$ is the (at most countable)
set of discontinuities
of the pointwise limit $\lambda$ of $\hs^{N-1}(\Gamma(\cdot))$ (which exists up to a further subsequence
by Helly's Theorem). Then we pass to the limit in \eqref{estabovec} obtaining
$$
\Es(t) \le \Es(0)+2\int_0^t \int_{\Om} \nabla u(\tau) \nabla \dot{g}(\tau) \,dx\,d\tau;
$$
moreover, following the proof of Theorem \ref{mainthm}, we have that
for all $t \in [0,1]$
$$
\Es(t) \le \liminf_n \Es_n(t) \le \limsup_n \Es_n(t)=\Es(0)
+2\int_0^t \int_{\Om} \nabla u(\tau) \nabla \dot{g}(\tau) \,dx\,d\tau,
$$
and taking into account \eqref{estbelowc} we get the convergence of the total
energy at any time. Since $\nabla u^n(t) \to \nabla u(t)$ strongly in
$L^2(\Om';\R^N)$ for every $t \in I_\infty$, we deduce that $\lambda=\hs^{N-1}(\Gamma(\cdot))$
on $I_\infty$, so that the convergence of the surface energy holds in $I_\infty$.
The extension to the continuity times for $\hs^{N-1}(\Gamma(\cdot))$ follows like in the final
part of the proof of Theorem \ref{mainthm}.

\bigskip
\bigskip
\centerline{ACKNOWLEDGMENTS}
\bigskip\noindent
The authors wish to thank Gianni Dal Maso for having proposed them the problem,
and Gianni Dal Maso and Gilles A. Francfort for many helpful and interesting discussions.

\end{document}